\documentclass[10pt,reqno,a4paper]{amsart}
\usepackage{euscript}
\usepackage[T1]{fontenc}
\usepackage{enumitem}
\usepackage{bm}
\usepackage{pgfplots}
\usepackage{float}
\usepackage{bbm}
\usepackage{amsthm}
\usepackage{amsmath}
\usepackage{amssymb}
\usepackage{subfigure}
\usepackage{tikz-cd}
\usepackage{graphicx,caption}
\usepackage{tikz}
\usetikzlibrary{arrows}
\usetikzlibrary{arrows.meta}
\usepackage{color}
\usepackage[bookmarksnumbered=true]{hyperref}
\hypersetup{
colorlinks = true,
linkcolor = black,
anchorcolor = black,
citecolor = black,
filecolor = black,
urlcolor = black
}
\newtheorem{theorem}{Theorem}
\newtheorem{proposition}[theorem]{Proposition}
\newtheorem{corollary}[theorem]{Corollary}
\newtheorem{lemma}{Lemma}
\theoremstyle{definition}
\newtheorem{remark}{Remark}
\newtheorem{example}{Example}

\renewcommand{\epsilon}{\varepsilon}
\def\R{\mathbb{R}}

\def\N{\mathbb{N}}
\def\Z{\mathbb{Z}}
\DeclareMathOperator{\Span}{span}
\DeclareMathOperator{\interior}{int}

\def\cM{\EuScript{M}}

\def\cU{\EuScript{U}}
\def\cW{\EuScript{W}}

\renewcommand{\epsilon}{\varepsilon}
\def\N{\mathbb{N}}
\def\R{\mathbb{R}}
\def\Z{\mathbb{Z}}
\def\cA{\EuScript{A}}
\def\cB{\EuScript{B}}
\def\cC{\EuScript{C}}
\def\cF{\EuScript{F}}

\def\cM{\EuScript{M}}
\def\cU{\EuScript{U}}

\def\cE{\EuScript{E}}

\def\cW{\EuScript{W}}

\newcommand{\M}{\mathcal{M}}

\begin{document}
	
\title[Ergodic full measures and entropy spectrum]{Ergodic full measures and entropy spectrum for countable observables on compact spaces}
	
\author[C. E. Holanda]{Carllos Eduardo Holanda}
\address{Institute of Mathematics, Universidade Federal de Alagoas, Macei\'o, 57072-970, Alagoas, Brazil}
\email{\tt c.eduarddo@gmail.com}
	
\keywords{Ergodic full measures, entropy spectrum, multifractal formalism, equilibrium states, countable observables}
	
\subjclass[2020]{28D20, 37D35.}
	
\begin{abstract}
We investigate conditions and obstructions for the existence of ergodic full measures in the context of countable observables on compact spaces. Extending the theory in a unified way, we are able to show a \emph{complete} picture via thermodynamic formalism for symbolic, uniformly hyperbolic and related systems. We prove conditional variational principles for the topological entropy, showing good approximations by ergodic full measures, and we give conditions where we precisely get ergodic full measures. In contrast with the case with a finite number of observables, new phenomena appear: we exhibit explicit examples of level sets with a null measure with respect to any invariant measure but having large entropy. On the other hand, based on induced ergodic systems in the sense of Kakutani-Rokhlin and heavily utilizing the formalism for countable Markov shifts, we go beyond the classical thermodynamic approach and show a  different way of building ergodic full measures. Additionally, we study the regularity and shape of the entropy spectrum and, in the end, we also give some applications of the main conditional variational principles. 
\end{abstract}

\maketitle
	

\section{Introduction}

Inspired by statistical physics, since its conception with D. Ruelle, R. Bowen and P. Walters (\cite{Rue73,Rue78,Bow75, Wal75}), thermodynamic formalism, developed inside ergodic theory, has a tendency to connect some different aspects of Mathematics, remarkably linking dimension theory, variational principles, equilibrium states and phase transitions, with many applications to number theory and surrounding areas. In this regard, this work is a further contribution to the ever expanding scope of the mathematical theory of thermodynamics and its various connections, here mainly involving multifractal analysis, together with aspects of rotation theory and also ergodic optimization.  

\subsection{The problem and the motivation}

Let $T\colon X \to X$ be a continuous map of a compact metric space $X$, and consider an \emph{observable} $\Phi \colon X \to \R^{\N}$ with continuous coordinate functions $\phi_j$. For each vector $\alpha = (\alpha_1,\alpha_2,\dots) \in \R^{\N}$, the \emph{level set} $C_{\alpha}$ is the set of points $x \in X$ such that the \emph{Birkhoff average} $(1/N)\sum_{k=0}^{N-1}\phi_j(T^{k}(x))$ converges exactly to $\alpha_j$ for all $j \ge 1$. A $T$-invariant measure $\nu$ is called a \emph{full measure for the level set $C_{\alpha}$}, or simply a \emph{full measure for $\alpha$}, if
\[
\int_{X}\phi_j\, d\nu = \alpha_{j} \quad \textrm{for all $j \ge 1$}, \qquad \nu(C_{\alpha}) = 1 \qquad \text{and} \qquad h_{\rm{top}}(C_{\alpha}) = h_{\nu}(T),
\]
where $h_{\rm{top}}(C_{\alpha}):= h_{\rm{top}}(T|C_{\alpha})$ is the \emph{Bowen topological entropy} of $T$ restricted to the level set, and $h_{\nu}(T)$ is the \emph{Kolmogorov--Sinai entropy} of $T$ with respect to $\nu$ (see Section \ref{INT} for the definitions). In this case we also say that $\alpha$ \emph{admits a full measure}. A full measure is thus like a conditional equilibrium state that sees the level set as a set of full measure and, at the same time, realizes its complexity. 

In this work we investigate the following single question:
\[
\textrm{\em{Are there (ergodic) full measures for countable observables on compact spaces?}}
\]
One motivating reason to ask the above question is the connection to the \emph{Katok conjecture} about intermediate entropy properties of diffeomorphisms (\cite{Kat80}), and some recent refinements via multifractal analysis (\cite{DHT24, BHHT25} and references therein). When $\Phi$ takes values in a finite-dimensional vector space, this question has a complete answer for many hyperbolic systems, obtained through the thermodynamic formalism (\cite{BS01,BS00,BSS02b,Cli13}), with some strong applications to number theory (\cite{BSS02a}). This is known today as the \emph{thermodynamic approach} to multifractal analysis (\cite{Cli14}). 

A crucial object in our work is the so called \emph{(generalized) rotation set} of $\Phi$ with respect to $T$, which is defined as
\[
\textrm{Rot}(\Phi) := \biggl\{\bigg(\int_{X}\phi_1\,d\mu, \int_{X}\phi_2\,d\mu,\cdots \bigg): \mu \in \cM(T)\biggr\} \subset \R^{\N},
\]
where $\cM(T)$ is the set of $T$-invariant probability measures. The finite-dimensional version of these sets, when $\Phi$ takes value in $\R^{d}$, are compact and convex sets that work as the domain of the \emph{entropy spectrum} $\alpha \mapsto h_{\rm{top}}(C_{\alpha})$. They are well known and have deep relations with thermodynamic formalism, ergodic optimization and rotation theory (see \cite{KW14,Jen19} and references therein).  

Starting in \cite{BS01}, and culminating in \cite{BSS02b}, Barreira, Saussol and Schmeling studied the entropy and dimension spectrum of $\R^{d}$-valued observables for uniformly hyperbolic and related systems and showed that, under a condition on uniqueness of equilibrium states, ergodic full measures always exist for vectors $\alpha$ in the \emph{interior} of the corresponding rotation set. Moreover, they are unique (\cite{BHHT25}) and are also unique equilibrium states for some explicit real-valued observables. Considering the same types of systems, but dropping the uniqueness assumption on the observable, Climenhaga \cite{Cli13,Cli14} substantially extended the machinery and showed that, for $\alpha$ in the interior of the rotation set, one still obtains conditional variational principles with sufficiently ``good'' approximations by ergodic full measures.

The setup with a compact space and at least countably many observables was already considered in the literature. Notably, Takens and Verbitskiy \cite{TV03}, Olsen \cite{Ols03}, and Pfister and Sullivan \cite{PS07} proved conditional variational principles for infinite-dimensional vector-valued observables, and Fan, Liao and Peyri\`ere \cite{FLP08} obtained them for Banach-valued observables, with some interesting applications. For setups where we do not necessarily have compact spaces, it is worth mentioning the more recent works \cite{FJLR15, JT21, Rus23}. However, many of these works developed the multifractal analysis via the \emph{orbit-gluing approach} (for saturated systems), and none of them considered the systematic study of (ergodic) full measures, or even measures supported on the level sets. For finite observables, the detailed information about full measures is a byproduct of the thermodynamic approach, and this is one strong motivating point of our treatment: to understand what survives, and what is lost, when the thermodynamic approach is scaled up to countably many observables.

\subsection{The results}

It turns out that searching for ergodic full measures for countable observables reveals new aspects of the theory, in particular new obstructions and nuances that are not present in the finite-dimensional framework. 

\subsubsection*{Conditional variational principles.} The starting point is that the entropy of a level set can be computed by finite-dimensional approximations. For a \emph{permissible system} $T\colon X \to X$, which includes many well known symbolic and uniformly hyperbolic maps on compact spaces (see the definition in Section \ref{INT}), and any $\alpha = (\alpha_1, \alpha_2, \dots) \in \textrm{Rot}(\Phi)$, we show in Theorem~\ref{ACVP1} that
\[
h_{\rm{top}}(C_{\alpha}) = \sup_{\mu \in \cM(\alpha)}h_{\mu}(T) = \lim_{N \to \infty}\sup_{\mu \in \cM_N(\alpha)}h_{\mu}(T),
\]
where $\cM_N(\alpha)$ is the set 
\[
\biggl\{\mu \in \cM(T): \biggl(\int_{X}\phi_1\,d\mu, \dots, \int_{X}\phi_N\,d\mu\biggr) = (\alpha_1,\dots,\alpha_N)\biggr\},
\]
and $\cM(\alpha) := \cap_{N \ge 1}\cM_N(\alpha)$. When $\alpha$ belongs to the \emph{pseudo-interior} of $\rm{Rot}(\Phi)$, that is, each vector $(\alpha_1,\dots,\alpha_N)$ is an interior point of the corresponding finite-dimensional rotation set in $\R^{N}$, the entropy of $C_\alpha$ is a limit of infima of topological pressures, and $C_\alpha$ is approximated, both in entropy and in measure, by ergodic equilibrium states of finite linear combinations of the coordinate functions $\phi_j$. This extends the main conditional variational principles in \cite{BSS02b,Cli13} to countable observables, and already brings more information than the orbit-gluing approach, now guaranteeing approximations by ergodic full measures and showing their complete characterization. Here, the main technical difficulty in the proofs is to control the ``thermodynamic" parameters uniformly with respect to the growing scale of approximation.

\subsubsection*{Existence and uniqueness: the regions $\cW$ and $\cA_{\rm{D}}$.} 
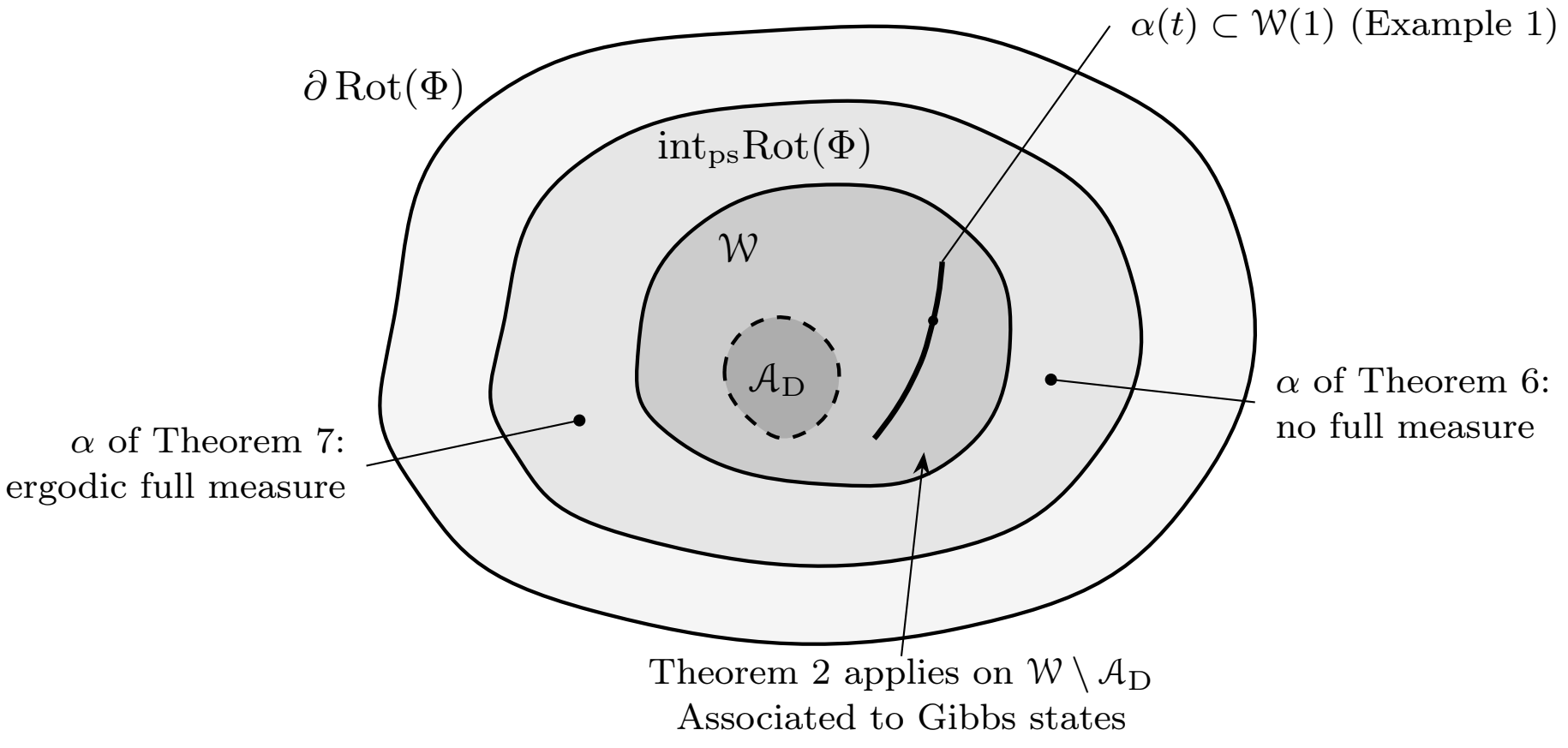
\begin{figure}[h!]
\centering
\begin{tikzpicture}[scale= 0.7, every node/.style={font=\small}]
\draw[thick, fill=black!4]
plot[smooth cycle, tension=0.8] coordinates
{(-4.6,0.2) (-3.6,2.6) (-0.4,3.4) (2.8,3.0) (4.7,1.2) (4.3,-1.6) (1.6,-3.2) (-2.4,-3.0) (-4.4,-1.6)};
\node at (-4.7,2.75) {$\partial\,\textrm{Rot}(\Phi)$};

\draw[thick, fill=black!10]
plot[smooth cycle, tension=0.8] coordinates
{(-3.4,0.1) (-2.6,1.9) (-0.2,2.6) (2.1,2.2) (3.5,0.8) (3.2,-1.2) (1.2,-2.4) (-1.9,-2.2) (-3.3,-1.2)};
\node at (-0.5,2.1) {$\textrm{int}_{\rm ps}\textrm{Rot}(\Phi)$};

\draw[thick, fill=black!20]
plot[smooth cycle, tension=0.8] coordinates
{(-1.9,-0.1) (-1.3,1.2) (0.3,1.7) (1.7,1.2) (2.2,0.0) (1.6,-1.3) (0.1,-1.6) (-1.5,-1.1)};
\node at (-0.8,1.0) {$\cW$};

\draw[very thick]
plot[smooth, tension=0.9] coordinates {(0.7,-1.1) (1.1,-0.5) (1.35,0.2) (1.45,0.85)};
\filldraw (1.35,0.2) circle (1.3pt);
\draw[thin] (1.45,0.85) -- (3.3,3.45);
\node[font=\footnotesize, right=1pt] at (3.3,3.45) {$\alpha(t)\subset\cW(1)$ (Example~1)};

\draw[thick, dashed, fill=black!32]
plot[smooth cycle, tension=0.8] coordinates
{(-0.6,-0.95) (-0.95,-0.35) (-0.55,0.2) (0.1,0.05) (0.3,-0.55) (-0.1,-1.05)};
\node at (-0.35,-0.45) {$\cA_{\rm D}$};

\filldraw (2.65,-0.45) circle (1.6pt);
\draw[thin] (2.65,-0.45) -- (4.9,-0.7);
\node[font=\footnotesize, align=left, right=1pt] at (4.9,-0.7)
{$\alpha$ of Theorem~6:\\ no full measure};

\filldraw (-2.55,-0.9) circle (1.6pt);
\draw[thin] (-2.55,-0.9) -- (-4.9,-1.4);
\node[font=\footnotesize, align=right, left=1pt] at (-4.9,-1.4)
{$\alpha$ of Theorem~7:\\ ergodic full measure};

\draw[-{Stealth[length=5pt]}] (1.0,-3.5) -- (1.25,-1.25);
\node[font=\footnotesize, align=center] at (1.0,-3.9)
{Theorem~2 applies on $\cW\setminus\cA_{\rm D}$\\ Associated to Gibbs states};
\end{tikzpicture}
\caption{The inclusions $\cA_{\rm D}\subset\cW\subset\textrm{int}_{\rm ps}\textrm{Rot}(\Phi)$ are proper in general (Theorems~6 and~7), while for finite-dimensional observables all three sets coincide with $\textrm{int}\,\textrm{Rot}(\Phi)$.}\label{BPIC}
\end{figure}

For permissible systems, passing from approximations to an actual ergodic full measure is where the countable setting departs from the finite one. In finite dimension,  because of the differentiability of the pressure, the interior of the rotation set is the exact optimal place to look for full measures (\cite{BS01,BSS02b}). With countably many observables, interior points alone are not enough, and we introduce two more layers inside the rotation sets (see Figure \ref{BPIC}). A vector $\alpha$ is \emph{uniformly attainable}, and we write $\alpha \in \cW$, if for each $N$ the finite projection $(\alpha_1,\dots,\alpha_N)$ is realized by the equilibrium state of some observable $\sum_{j \le N}q_j\phi_j$ with the vector $q = (q_1,\dots,q_N)$ satisfying $\|q\|_{\infty} \le R$, with the bound $R$ being independent of $N$. A vector is \emph{asymptotically deep}, and we write $\alpha \in \cA_{\rm{D}}$, if its finite projections stay at a positive distance from the boundaries of the finite-dimensional rotation sets in $\R^{N}$, uniformly in the dimension $N \ge 1$. The notion of $\cA_{\rm{D}}$ is more geometric in nature, and we show that $\cA_{\rm{D}} \subseteq \cW$ (Proposition~\ref{PP}). Our main existence result via thermodynamic formalism, Theorem~\ref{ACVP2}, states that if $\alpha \in \cW$ and the observables $\sum_{j\ge1}\gamma_j(\phi_j - \alpha_j)$ with $\|\gamma\|_\infty \le R$ have uniformly small tails and unique equilibrium states, then $\alpha$ admits a \emph{unique} ergodic full measure, which is the unique equilibrium state for 
\[
\varphi:= \sum_{j \ge 1}q^{*}_j(\phi_j - \alpha_j) \quad \textrm{for some $q^{*} = (q^{*}_1, q^{*}_2, \dots)\in \R^{\N}$ with $\|q^{*}\|_{\infty} \le R$},
\]
and $h_{\rm top}(C_\alpha)$ is equal to the topological pressure of $\varphi$. The proof is a compactness argument on the parameters combined with the stability of equilibrium states. Perhaps even counterintuitively we also show that the above tail condition, appearing as an hypothesis in Theorem \ref{ACVP2}, fails on $\cA_{\rm D}$ (Proposition \ref{OSC}). This means that the mechanism of Theorem \ref{ACVP2} can only work for vectors $\alpha \in \cW \setminus \cA_{\rm D}$, the ones that are, a priori, not very deep in the rotation set.  

\subsubsection*{Outside $\cW$: non-existence.} Our first main example, Theorem \ref{NEFM}, shows that the thermodynamic method cannot be extended beyond $\cW$ by any refinement of its hypotheses. For the full shift on two symbols, we construct an observable $\Phi$ with locally constant coordinates satisfying all the regularity assumptions of Theorem \ref{ACVP2}, and a vector $\alpha$ in the pseudo-interior of $\textrm{Rot}(\Phi)$ but outside $\cW$, with the following properties: the level set $C_\alpha$ has topological entropy close to the absolute maximum $h_{\rm top}(\sigma)$, and yet $\mu(C_\alpha) = 0$ for \emph{every} invariant measure $\mu$. In particular, $\alpha$ does not admit any full measure, ergodic or not. The construction, inspired in parts by some techniques from \cite{Bou01}, is purely infinite-dimensional. In fact, the coordinate functions of the observable are chosen  level by level, so that they separate invariant measures and no nontrivial finite combination of them is cohomologous to a constant. For a finite number of observables this situation would never happen, as non-ergodic measures realizing an interior point is always well approximated by ergodic ones, and with positive entropy (\cite{Cli13}). 

\subsubsection*{Outside $\cW$: existence by a different mechanism.} Going in the opposite direction, our second main example, Theorem \ref{NTW}, shows that $\cW$ is not necessary for the existence of ergodic full measures either. However, outside $\cW$ they must be produced by a different mechanism other than that of Theorem \ref{ACVP2}. Again for the full shift on two symbols, we consider the observable with coordinate functions $\mathbf 1_{A_j}$, where $A_j$ is the cylinder generated by the word $0^{j-1}1$, and the vector $\alpha$ is given by $\alpha_j = c\,4^{-2^{j}}$ (with the constant $c \in \R$ chosen to satisfy $\sum_{j \ge 1}\alpha_j = 1$). We show that $\alpha$ is in the pseudo-interior of the rotation set but not in $\cW$. Nevertheless, we are able to show that $\alpha$ admits a unique ergodic full measure, with $0 < h_{\rm top}(C_\alpha) < \log 2$. This measure is built by invoking \emph{Kakutani-Rokhlin towers}, where the induced dynamics behaves as a full shift on countably many symbols and is linked to a prescribed probability vector, which gives the desired weight for the cylinders $A_j$. Because of that, the construction borrows many tools from Sarig's formalism for countable topological Markov shifts (\cite{Sar99,Sar01}), a non-compact framework, even though we are working on a compact space. The entropy of the level set $C_{\alpha}$ is computed with a conditional entropy argument on the countable partition $\{A_1, A_2,\dots\}$, and it is attained exactly by the prescribed induced measure. In spite of the vector $\alpha$ be out of reach of the classical thermodynamic method, in the sense that the parameters of the approximating equilibrium states diverge, it is still realized by a measure which, on the induced system, is itself an equilibrium state. 

We believe that the aforementioned mechanism (induced systems together with prescribed probability vectors on $\N$) is interesting on its own, and might be applicable in other contexts. Figure \ref{BPIC} summarizes the resulting relations between the three layers inside the rotation set
\[
\cA_{\rm D} \subset \cW \subset \interior_{\rm ps}\textrm{Rot}(\Phi),
\]
with both inclusions proper in general. Theorem \ref{ACVP2} holds in $\cW \setminus \cA_{\rm D}$, and the two main examples work in the region $\interior_{\rm ps}\textrm{Rot}(\Phi) \setminus \cW$, but with opposite behaviors. In Section \ref{NRC} we show a case where we can explicitly calculate $\cW(R)$ for a suitable parameter $R > 0$. For finite observables, all three sets coincide with the interior of the finite-dimensional rotation set.

\subsubsection*{Regularity of the spectrum.} For the full shift and the observable tracking the frequencies of all possible cylinders, Theorem \ref{REG} shows that the spectrum $\alpha \mapsto h_{\rm top}(C_\alpha)$ is affine and upper semicontinuous on the rotation set, vanishes on a dense $G_\delta$ set, where it is continuous, and is discontinuous at every point of positive entropy, which also form a dense set. Only upper semicontinuity and continuity on convex paths (Example \ref{CPTH}) survive from the analytic spectrum obtained in the setup with finite observables.  The reason is again separation of invariant measures, forcing the integration map to become an affine homeomorphism between $\cM(\sigma)$ and the rotation set. In this case, the spectrum is actually reduced to the metric entropy map transported by it.

\subsubsection*{Applications.} We apply the conditional variational principles (Theorems \ref{ACVP1} and \ref{ACVP2}) to the symbolic distribution of infinite words, that is, to level sets of prescribed frequencies of countably many cylinders, complementing the results of \cite{FLP08}. More specifically, for a parameter in the pseudo-interior of the rotation set, the maximal entropy over the level set is approximated by ergodic equilibrium states (Corollary~\ref{BEF}), and for suitably normalized frequencies, uniformly attainable vectors admit unique ergodic full measures (Section \ref{NRC}). Finally, as a byproduct of the thermodynamic method, we show how to characterize the measures that realize uniformly attainable vectors with an infinite number of linear constraints in the presence of a reference potential. They are Gibbs states for an observable with countably many parameters (Corollary \ref{GIBBS}), then extending Theorem G in \cite{GKLM18}. \\

For reference, putting now everything together, we can see the complete landscape that captures the new phenomena with countable observables. 

\subsection*{Structure of the paper} Section \ref{INT} collects the needed background on thermodynamic formalism and Bowen's entropy. Section \ref{MFA} contains the main results: the conditional variational principles and the existence theorem in Section \ref{CVPS}, the comparison between the geometric and the thermodynamic conditions in Section \ref{GCUA}, the two examples outside $\cW$ in Sections \ref{NW} and \ref{OTW}, and the regularity of the spectrum in Section \ref{RES}. Section \ref{ACV} contains the applications. 

\subsection*{Acknowledgments.} This work started in a visit to Fudan University in Shanghai, China, in December 2023. In this regard, the author is grateful to Xueting Tian and Xiaobo Hou for some fruitful conversations about the topic, and possible applications of ergodic full measures to intermediate entropy problems. I would also like to thank Al\'icia Mayumi for her daily presence during the preparation of this manuscript. The author is supported by Conselho Nacional de Desenvolvimento Científico e Tecnológico (CNPq) project 446515/2024-8.

\subsection*{A.I. disclosure} Large Language Models were used to proofread parts of the manuscript, to assist in some computations and in literature search, and also to generate the required figures. All the main ideas, structure of the work, and proofs are the responsibility of the author only. 
	
\section{Preliminaries on thermodynamic formalism on compact spaces}\label{INT}

Let $X$ be a compact metric space, $T\colon X\to X$ a continuous map, and denote $\mathcal{M}(T)$ as the set of $T-$invariant
probability measures on $X$, which we always regard with the weak$^{*}$ topology. We write $\mathcal{M}_{\rm{erg}}(T)$ for the subset of ergodic measures. For any function $\varphi: X \to \R$, we denote the Birkhoff sums by
\[
S_{N}\varphi:=\sum_{k=0}^{N-1}\varphi\circ T^{k}.
\]
Throughout the text, as in the related literature, we might also say that $\varphi$ is an \emph{observable} or a \emph{potential}. \\  

\subsection*{Kolmogorov-Sinai (metric) entropy.} Let $\mu\in\mathcal{M}(T)$. Given a finite partition $\xi=\{C_{1},\dots,C_{m}\}$ of $X$ into Borel sets, we define its entropy as 
\[
H_{\mu}(\xi):=-\sum_{i=1}^{m}\mu(C_{i})\log\mu(C_{i}).
\]
Given two partitions $\xi$ and $\eta$, we denote by $\xi\vee\eta$ the partition formed by the intersections
$C\cap D$ with $C\in\xi$ and $D\in\eta$, and we write
\[
\xi^{\,n}:=\bigvee_{i=0}^{n-1}T^{-i}\xi = \xi \vee T^{-1}\xi \vee \cdots \vee T^{n-1}\xi.
\]
The sequence $\bigl(H_{\mu}(\xi^{\,n})\bigr)_{n\ge1}$ is subadditive, so the limit below exists and coincides with the infimum, that is, 
\[
h_{\mu}(T,\xi):=\lim_{n\to\infty}\frac1n H_{\mu}(\xi^{\,n})
=\inf_{n\ge1}\frac1nH_{\mu}(\xi^{\,n}).
\]
The \emph{Kolmogorov-Sinai entropy}, also called \emph{metric entropy} of $(T,\mu)$, is defined as
\[
h_{\mu}(T):=\sup\bigl\{h_{\mu}(T,\xi):\ \xi\ \text{a finite Borel partition of }X\bigr\}.
\]
The map $\mu\mapsto h_{\mu}(T)$ is affine on $\mathcal{M}(T)$. 

We say that a partition $\xi$ is generating for $T$ when $\bigvee_{n\ge0}T^{-n}\xi$ coincides, modulo $\mu-$null sets, with the Borel $\sigma-$algebra of $X$. When $\xi$ is a finite generating partition, \emph{Kolmogorov-Sinai} theorem gives $h_{\mu}(T)=h_{\mu}(T,\xi)$.\\

\subsection*{Topological pressure and topological entropy.} We now measure the complexity of $T$ itself, with a weight given by a
potential. For $n\ge1$ we consider the Bowen metric
\[
d_{n}(x,y):=\max_{0\le i<n}d\bigl(T^{i}x,T^{i}y\bigr),
\]
and, for $\varepsilon>0$, the Bowen ball
\[
B_{n}(x,\varepsilon):=\{y\in X:\ d_{n}(x,y)<\varepsilon\},
\]
which is the set of points whose orbit stays $\varepsilon-$close to the orbit
of $x$ during $n$ iterates. A set $E\subseteq X$ is said to be
$(n,\varepsilon)-$separated when $d_{n}(x,y)>\varepsilon$ for every pair of
distinct points $x,y\in E$. Since $X$ is compact, such a set is necessarily
finite. Given a continuous potential $\varphi: X \to \R$, we set
\[
Z_{n}(\varphi,\varepsilon):=\sup\Bigl\{\sum_{x\in E}e^{S_{n}\varphi(x)}:\
E\ \text{is an }(n,\varepsilon)-\text{separated subset of }X\Bigr\},
\]
and we define the \emph{topological pressure} of $\varphi$ as
\[
P_{\mathrm{top}}(\varphi)
:=\lim_{\varepsilon\to0}\ \limsup_{n\to\infty}
\frac1n\log Z_{n}(\varphi,\varepsilon).
\]
The limit in $\varepsilon$ exists because the quantity inside is non-decreasing as $\varepsilon$ decreases. The \emph{topological entropy} of $T$ is the particular case with $\varphi \equiv 0$
\[
h_{\mathrm{top}}(T):=P_{\mathrm{top}}(0).
\]
It is easy to see from the definition that $P_{\mathrm{top}}$ is convex and Lipschitz with respect to the supremum norm $\|\cdot\|_{\infty}$, that is, for each pair of continuous potentials $\varphi, \psi$ we have
\[
\bigl|P_{\mathrm{top}}(\varphi)-P_{\mathrm{top}}(\psi)\bigr|
\le\|\varphi-\psi\|_{\infty}.
\]
Moreover, $P_{\mathrm{top}}$ is finite whenever
$h_{\mathrm{top}}(T)<\infty$.\\

\subsection*{Variational principle and equilibrium states.} The metric entropy and the topological pressure are connected via the well known \emph{variational principle}
\[
P_{\mathrm{top}}(\varphi)
=\sup_{\mu\in\mathcal{M}(T)}\Bigl\{h_{\mu}(T)+\int_{X}\varphi\,d\mu\Bigr\}
=\sup_{\mu\in\mathcal{M}_{\rm{erg}}(T)}\Bigl\{h_{\mu}(T)+\int_{X}\varphi\,d\mu\Bigr\}.
\]
In particular $h_{\mathrm{top}}(T)=\sup_{\mu\in\mathcal{M}(T)}h_{\mu}(T)$. 

We say that $\mu\in\mathcal{M}(T)$ is an \emph{equilibrium state} for $\varphi$ when it realizes the supremum in the variational principle above, that is, when
\[
h_{\mu}(T)+\int_{X}\varphi\,d\mu=P_{\mathrm{top}}(\varphi).
\]
For more information on the thermodynamic formalism on compact spaces, see for example, \cite{VO16} and references therein. \\

\subsection*{Bowen topological entropy.} The sets whose complexity we want to measure, such as the level sets $C_{\alpha}$, are in general neither compact nor $T-$invariant, so the definition above does not apply to them directly. We use instead the Carathéodory type construction of Bowen \cite{Bow73, Bow75}. 

Let $Z\subseteq X$ be an arbitrary subset. Fix $\varepsilon>0$, $s\ge0$ and $N\in\mathbb{N}$, and set
\[
M(Z,s,N,\varepsilon):=\inf\ \sum_{i \ge 1}e^{-s\,n_{i}},
\]
where the infimum is taken over all finite or countable families
$\{B_{n_{i}}(x_{i},\varepsilon)\}_{i}$ of Bowen balls, with $x_{i}\in X$ and $n_{i}\ge N$, whose union covers $Z$. The quantity $M(Z,s,N,\varepsilon)$ is
non-decreasing in $N$, so that the limit
\[
m(Z,s,\varepsilon):=\lim_{N\to\infty}M(Z,s,N,\varepsilon)
\]
exists in $[0,+\infty]$. As a function of $s$ it jumps from $+\infty$ to $0$
at a single critical value, and this value is what we call
\[
h_{\mathrm{top}}(T|Z,\varepsilon)
:=\inf\{s\ge0:\ m(Z,s,\varepsilon)=0\}
=\sup\{s\ge0:\ m(Z,s,\varepsilon)=+\infty\}.
\]
This quantity is non-decreasing as $\varepsilon$ decreases, and therefore the \emph{Bowen topological entropy} of $T$ on $Z$ is defined as
\[
h_{\mathrm{top}}(T|Z):=\lim_{\varepsilon\to0}h_{\mathrm{top}}(T|Z,\varepsilon).
\]
It is monotone, in the sense that $h_{\mathrm{top}}(T|Z_{1})\le
h_{\mathrm{top}}(T|Z_{2})$ whenever $Z_{1}\subseteq Z_{2}$, and countably
stable, that is
$h_{\mathrm{top}}\bigl(T|\bigcup_{i}Z_{i}\bigr)=\sup_{i}h_{\mathrm{top}}(T|Z_{i})$.
More specifically, when $Z$ is compact and $T-$invariant it coincides with
the topological entropy of the restriction $T|_{Z}$, and in particular
$h_{\mathrm{top}}(T|X)=h_{\mathrm{top}}(T)$. 

\section{Multifractal formalism via finite-dimensional approximations}\label{MFA}
	
In this section we start with an approach of approximating infinite dimensional systems with finite dimensional ones and invoking known results of multifractal formalism via thermodynamics, that is, concerning the connections with the topological pressure and existence (or good approximations) of ergodic measures supported on the level sets. However, as we shall see, this cannot be done directly and the main difficulty is to make sure we can approach infinite dimensional systems in the right way, transporting the information about entropies and dimensions with a ``good'' dependency on the approximation scale.  
	
\subsection{Conditional variational principles and ergodic full measures}\label{CVPS}
We start fixing some notations. $X$ is a compact metric space, $T\colon X \to X$ is a continuous map, and $\Phi$ is a $\R^{\N}-$valued observable with continuous coordinate functions $\phi_j$. Given a vector $\alpha = (\alpha_1, \alpha_2, \dots) \in \rm{Rot}(\Phi)$ and $N \ge 1$, we set
\[
\cM_N(\alpha) := \biggl\{\mu \in \cM(T): \biggl(\int_{X}\phi_1\,d\mu, \dots, \int_{X}\phi_N\,d\mu\biggr) = (\alpha_1,\dots,\alpha_N)\biggr\},
\]
and define $\cM(\alpha) = \cap_{N \ge 1}\cM_N(\alpha)$. 
	
\begin{lemma}\label{USC} Assume that the metric entropy map $\cM(T) \ni \mu \mapsto h_{\mu}(T)$ is upper semicontinuous. Then, for any $\alpha \in \rm{Rot}(\Phi)$ we have
\[
\sup_{\mu \in \cM(\alpha)}h_{\mu}(T) = \lim_{N \to \infty}\sup_{\mu \in \cM_N(\alpha)}h_{\mu}(T). 
\]
\end{lemma}
	
\begin{proof}
Let $\alpha \in \rm{Rot}(\Phi)$. Since $\cM_{N+1}(\alpha) \subset \cM_{N}(\alpha)$ for every $N \ge 1$, we get $\sup_{\mu \in \cM_{N+1}(\alpha)}h_{\mu}(T) \le \sup_{\mu \in \cM_{N}(\alpha)}h_{\mu}(T)$. This means that the sequence $A_N:= \sup_{\mu \in \cM_{N}(\alpha)}h_{\mu}(T)$ is non-increasing. Because $\cM(\alpha) \subset \cM_{N}(\alpha)$ for every $N$, we have $A_N \ge \sup_{\mu \in \cM(\alpha)}h_{\mu}(T)$ for all $N \ge 1$. Consequently, the limit $A:= \lim_{N \to \infty}A_N$ exists, it is equal to $\inf_{N \ge 1}A_N$, and satisfies $A \ge \sup_{\mu \in \cM(\alpha)}h_{\mu}(T)$. $A$ might be infinite.  By contradiction, assume that $A$ is strictly larger than $\sup_{\mu \in \cM(\alpha)}h_{\mu}(T)$. Now pick some $\delta$ with $0 < \delta < A - \sup_{\mu \in \cM(\alpha)}h_{\mu}(T)$. By the definition of the sequence $(A_N)_{N \ge 1}$, for each $N$ there exists $\eta_N \in \cM_{N}(\alpha)$ such that $h_{\eta_{N}}(T) > A_N - 1/N$, which readily gives 
\begin{equation}\label{INF}
\liminf_{N \to \infty}h_{\eta_N}(T) \ge A.
\end{equation}
On the other hand, the sequence $(\eta_N)_{N \ge 1} \subset \cM_1(\alpha)$. The compactness of $\cM_1(\alpha)$ guarantees the existence of a convergent subsequence $(\eta_{N_k})_{k \ge 1}$ with limit $\eta_* \in \cM_1(\alpha)$. Observe that for each arbitrary $M$ if $N_k \ge M$ then $\eta_{N_k} \in \cM_{N_k}(\alpha) \subset \cM_M(\alpha)$, implying that $\eta_*$ belongs to the closure of $\cM_{M}(\alpha)$, which is $\cM_M(\alpha)$ itself. Hence, $\eta_{*} \in \cap_{M \ge 1}\cM_M(\alpha) = \cM(\alpha)$. Now invoking the upper semicontinuity of the entropy map, we get 
\[
\limsup_{k \to \infty}h_{\eta_{N_k}}(T) \le h_{\eta_{*}}(T). 
\]
Together with \eqref{INF}, this gives
\[
h_{\eta_*}(T) \ge \limsup_{k \to \infty}h_{\eta_{N_k}}(T)\ge A.  
\]
However, $\eta_* \in \cM(\alpha)$ and we must have $A \le h_{\eta_*}(T) \le \sup_{\mu \in \cM(\alpha)}h_{\mu}(T)$, which is a contradiction. 
\end{proof}
	
Lemma \ref{USC} paves the way for the development of conditional variational principles via thermodynamic formalism using approximations by finite-dimensional  systems. As mentioned earlier, this unlocks many properties that were not possible to obtain by the orbit-gluing approach. 
	
For a given $\nu \in \cM(T)$, the set of \emph{$\nu-$generic points} is defined as
\[
G_{\nu}:=\biggl\{x \in X : \frac{1}{N}\sum_{k=0}^{N-1}\delta_{T^{k}(x)} \to \nu \biggr\},
\]
where $\delta_x$ is the \emph{Dirac measure} concentrated on $x$, and the convergence is in the weak$^* - $ topology. A dynamical system $T \colon X  \to X$ is said to be \emph{saturated} if for any $\mu \in \cM(T)$ we have $h_{\rm{top}}(T|G_{\mu}) = h_{\mu}(T)$. 
	
Now we are able to introduce the dynamical systems we are mainly interested in this work. We say that $(X,T)$ is a \emph{permissible system} if:
	
\begin{itemize}
\item $(X,T)$ is saturated.
		
\item $T|X$ has finite topological entropy. 
		
\item The metric entropy map $\mu \mapsto h_{\mu}(T)$ is upper semicontinuous. 
		
\item There is a dense set $\cU \subset C(X,\R)$ such that each $\xi \in \cU$ has a unique equilibrium state with respect to $T$ on $X$. 
\end{itemize}
	
\begin{remark}\label{PSYST}
Permissible systems include several types of well known dynamical systems, such as topologically mixing subshifts of finite type, topologically mixing $C^1$ maps that are uniformly expanding, topologically mixing uniformly hyperbolic $C^1$ maps and $\beta-$shifts satisfying the specification property (\cite{FLP08,KH95}). In all of these setups, one can take $\cU$ as the space of H\"older (Lipschitz) continuous functions, since these types of observables have unique equilibrium states, and are dense in the space of continuous functions on compact metric spaces.
\end{remark}
	
Let $\kappa$ be the metric on $\R^{\N}$ given by $\kappa(x,y) = \sum_{j \ge 1}2^{-j}\min\{1,|x_j - y_j|\}$. For any $\delta > 0$ we consider the set 
\[
C^{\delta}_{\alpha} := \biggl\{x \in X: \limsup_{N \to \infty}\kappa\bigg(\frac{1}{N}S_N\Phi(x),\alpha\bigg)\le \delta \biggr\} \supseteq C_{\alpha},
\] 
where $S_N\Phi(x) = \big(S_N\phi_1(x),S_N\phi_2(x),\dots\big)$. Note that $\cap_{\delta > 0}C^{\delta}_{\alpha} = C_{\alpha}$.  
    
Also recall that a function $f \colon X \to \R$ is \emph{cohomologous} to a function $g \colon X \to \R$ (with respect to~$T$) if there exists a measurable bounded function $q\colon X \to \R$ such~that
\[
f(x) = g(x) + q(T(x)) - q(x)\quad\text{for} \ x \in X.
\]
To simplify the notation, we denote $\phi^{(N)}$ as the projected $\R^{N}-$valued observable with coordinate functions $\phi_1,\dots,\phi_N$, and $v^{(N)} := (v_1,\dots,v_N)$, for any vector $v = (v_1,v_2,\dots) \in \R^{\N}$. 
     
For convenience, we also write
\[
\interior_{\rm{ps}}\textrm{Rot}(\Phi) := \big\{\alpha \in \textrm{Rot}(\Phi): \alpha^{(N)} \in \interior \textrm{Rot}(\phi^{(N)}) \textrm{ for all $N \ge 1$}\big\}
\]
and call it the \emph{pseudo-interior} of the rotation set. Ocasionally, we also may use the notation $\int_{X}\Phi\, d\mu = \big(\int_{X}\phi_1\,d\mu, \int_{X}\phi_2\,d\mu, \dots\big)$.
	
The following result gives conditions for ``good'' approximations by ergodic full measures. 
	
\begin{theorem} \label{ACVP1} Let $(X,T)$ be a permissible system and $\Phi\colon X \to\R^{\N}$ an observable with continuous coordinate functions $\phi_j$ that are continuous.  
\begin{enumerate}
\item For $\alpha \in \rm{Rot}(\Phi)$, we have
\[
h_{\rm{top}}(C_{\alpha}) = \lim_{N \to \infty}\sup_{\mu \in \cM_N(\alpha)}h_{\mu}(T) = \sup_{\mu \in \cM(\alpha)}h_{\mu}(T). 
\]
			
\item For $\alpha \in \interior_{\rm{ps}}\rm{Rot}(\Phi)$, we have
\[
h_{\rm{top}}(C_{\alpha}) = \lim_{N \to \infty}\inf_{q^{(N)} \in \R^N}P_{\rm{top}}\bigg(\sum_{j\le N}q_j (\phi_j-\alpha_j)\bigg),
\]
and for a given $\delta > 0$ there exists an ergodic measure $\nu_{\alpha} \in \cM(T)$ satisfying $\nu_{\alpha}(C^{\delta}_{\alpha}) = 1$, $\kappa(\int_X \Phi\,d\nu_{\alpha},\alpha) < \delta$ and such that
\[
\big|h_{\rm{top}}(C_{\alpha}) - h_{\nu_{\alpha}}(T)\big| < \delta \textrm{ and } \lim_{\delta \to 0}h_{\rm{top}}(C^{\delta}_{\alpha}) = h_{\rm{top}}(C_{\alpha}). 
\]
\end{enumerate}
\end{theorem}
	
\begin{proof}
The second equality of item 1 is a direct application of Lemma \ref{USC}. Since the system is in particular saturated, the first equality is Proposition 7.1 in \cite{PS07}, which is also the main result in \cite{TV03}.   

The relation with the topological pressure in item 2 follows from item 1 together with Theorem A in \cite{Cli13}. 

We now start proving the second part involving approximations by ergodic full measures. Fix $\delta>0$. Choose $N_{1}$ with $2^{-N_{1}}<\delta$ and, by item 1, choose $N_{2}$ with $|h_{\rm{top}}(C_{\alpha})- h_{\rm{top}}(C^{(N)}_{\alpha})|<\delta/2$ for all $N\ge N_{2}$. Take $N^* :=\max(N_{1},N_{2})$. By assumption $\alpha \in \interior_{\rm{ps}}\rm{Rot}(\Phi)$, which means $(\alpha_1,\dots,\alpha_N) \in \interior \rm{Rot}$$(\phi^{(N)})$ for each $N \ge 1$. In particular, for each $N \ge N^*$ Theorem 3.3 in \cite{Cli13} guarantees the existence of an ergodic measure $\nu_{\alpha}\in\mathcal M_{N}(\alpha)$, the unique equilibrium state of $\sum_{j\le N}q_{j}(\phi_{j}-\alpha_{j})$ for some $q\in\R^{N}$, with $\nu_{\alpha}(C^{(N)}_{\alpha})=1$ and satisfying $|h_{\rm{top}}(C^{(N)}_{\alpha})-h_{\nu_{\alpha}}(T)|<\delta/2$, where $C^{(N)}_{\alpha}$ is taken as the level set of $\phi^{(N)}$ at the vector $\alpha^{(N)}$. Thus, the first entropy estimate of item 2 follows now by a triangle inequality.

Write $b_{\nu_{\alpha}}:=\int_X\Phi\,d\nu_{\alpha}$, which is defined here coordinate-wise. Since $\nu_{\alpha}\in\mathcal M_{N}(\alpha)$, the first $N$ summands of $\kappa(b_{\nu_{\alpha}},\alpha)$ vanish, while each remaining summand is at most $2^{-j}$ because $\min(1,\cdot)\le1$. Hence
\begin{equation}\label{eq:tail}
\kappa(b_{\nu_{\alpha}},\alpha)\;\le\;\sum_{j>N}2^{-j}\;=\;2^{-N}\;<\;\delta.
\end{equation}
Since $\nu_{\alpha}$ is ergodic, Birkhoff's ergodic theorem applied to each $\phi_{j}$, together with a countable union of full sets, gives for $\nu_{\alpha}-$almost every $x \in X$
\[	
\frac1NS_{N}\phi_{j}(x)\longrightarrow\int_X\phi_{j}\,d\nu_{\alpha}
\qquad\text{for every }j\ge1,
\]
that is, $\frac{1}{N}S_{N}\Phi(x)\to b_{\nu_{\alpha}}$ in the product topology of
$\R^{\N}$. As $\kappa(\cdot,\alpha)$ is continuous for that topology, $\kappa\bigl(\frac1NS_{N}\Phi(x),\alpha\bigr)\to\kappa(b_{\nu_{\alpha}},\alpha)<\delta$ by \eqref{eq:tail}, so $\nu_{\alpha}(C^{\delta}_{\alpha})=1$. 

Now define $A_{\delta} = \{\mu \in \cM(T): \kappa\big(\int_{X}\Phi\,d\mu,\alpha\big)\le \delta\}$. The upper semicontinuity of the metric entropy allows us to invoke Lemma \ref{USC}, with $A_{\delta}$ in place of $\cM_{N}(\alpha)$, together with item 1, to conclude that
\begin{equation}\label{LBM}
\sup_{\mu \in A_{\delta}}h_{\mu}(T) \to \sup_{\mu \in \cM(\alpha)}h_{\mu}(T) = h_{\rm{top}}(C_{\alpha}) \quad \textrm{when $\delta \to 0$.}
\end{equation}
We also note that 
\[
C^{\delta}_{\alpha} =\{x \in X: V(x) \subseteq A_{\delta}\},
\]
where $V(x)$ is the set of limit points of the empirical measures $\frac{1}{N}\sum_{k=0}^{N-1}\delta_{T^{k}(x)}$. It follows from Theorem 2 in \cite{Bow73} that 
\[
h_{\rm{top}}(C^{\delta}_{\alpha}) = h_{\rm{top}}(\{x \in X: V(x) \subseteq A_{\delta}\})\le \sup_{\mu \in A_{\delta}}h_{\mu}(T).
\]
The lower bound clearly follows from saturation. Then, together with \eqref{LBM}, we finally get $h_{\rm{top}}(C^{\delta}_{\alpha}) \to h_{\rm{top}}(C_{\alpha})$ as $\delta \to 0$, and the result is proved. 
\end{proof}
The following subsets are based on the thermodynamic method via regularity of the topological pressure (\cite{BS01,BSS02b}), and are crucial for understanding where to find full measures in our infinite-dimensional framework. Given $R > 0$, $N\ge1$ and $q^{(N)} \in\R^{N}$, we let $\nu_{q}$ denote the unique equilibrium state of the real-valued function $\sum_{j\le N}q_{j}\phi_{j}$ (when it exists) and set
\begin{align*}
\cW_{N}(R)&:=\bigg\{\bigg(\int\phi_{1}\,d\nu_{q},\dots,\int\phi_{N}\,d\nu_{q}\bigg):\ \|q^{(N)}\|_{\infty}\le R\bigg\}\subseteq \textrm{Rot}(\phi^{(N)}),\\ \cW(R)&:=\bigl\{\alpha \in \R^{\N}:\ \alpha^{(N)} \in\mathcal \cW_{N}(R)\ \text{for all }N\ge1\bigr\}.
\end{align*}
A vector $\alpha \in \cW:= \bigcup_{R > 0}\cW(R)$ is said to be \emph{uniformly attainable}.  
\begin{remark}
We note that $\alpha\in \cW(R)$ guarantees, for each $N$ separately, the existence of a vector $q^{(N)}\in\R^{N}$ with $\|q^{(N)}\|_{\infty}\le R$ whose equilibrium state has $\int_X\phi_{j}\,d\nu_{q^{(N)}}=\alpha_{j}$ for $j\le N$. The bound $R$ is uniform in $N$, but the vectors themselves need not be compatible in the following sense: nothing forces $q^{(N+1)}_{j}=q^{(N)}_{j}$, since the conditions $\alpha^{(N)}\in \cW_{N}(R)$ are imposed independently at each level.
\end{remark}
	
The following result gives conditions for the existence of ergodic full measures.
 
\begin{theorem}\label{ACVP2} Under the hypotheses of Theorem \ref{ACVP1}, let $R > 0$, set $Q_{R}:= \{\gamma \in \R^{\N}: \sup_{j}|\gamma_j|\le R\}$, and assume that:
\begin{itemize}
\item[\textup{(A1)}] $\alpha \in \cW(R)$. 
\item[\textup{(A2)}] The uniform tail sequence 
\[
\varepsilon_{N}:= \sup_{\gamma \in Q_{R}}\bigg\|\sum_{j > N}\gamma_j(\phi_j-\alpha_j)\bigg\|_{\infty} \longrightarrow 0 \quad \textrm{when $N \to \infty$}.
\]
\item[\textup{(A3)}] Every $\psi_{\gamma}:= \sum_{j \ge 1}\gamma_j(\phi_j - \alpha_j)$ with $\gamma \in Q_{R}$ has a unique equilibrium state\footnote{These observables are not required to be necessarily in the set $\cU$ as in the definition of permissible systems, but they still have unique equilibrium states.}. 
\end{itemize}
Then $\alpha$ admits a unique ergodic full measure. Moreover, 
\[
 h_{\rm{top}}(C_{\alpha}) = P_{\rm{top}}\bigg(\sum_{j\ge1}q^*_j (\phi_j-\alpha_j)\bigg) \quad \textrm{for some $q^{*} \in Q_R$}. 
\] 
\end{theorem}
\begin{proof}
Since $\alpha \in \cW(R)$, for each $N \ge 1$ there exists a vector $q^{(N)} = (q^{(N)}_1,\dots,q^{(N)}_N) \in \R^{N}$ with $\|q^{(N)}\|_{\infty} \le R$ and such that $\int_{X}\phi_j\,d\nu_{N} = \alpha_j$ for each $j \le N$, where $\nu_{N}:=\nu_{q^{(N)}}$ is, in particular, the unique equilibrium state of the function
\[
\sum_{j \le N}q_j^{(N)}(\phi_j-\alpha_j) - h_{\rm{top}}(C^{(N)}_{\alpha}).
\]
One can check that $\nu_{N}$ is ergodic, belongs to $\cM_N(\alpha)$, satisfies $\nu_{N}(C^{(N)}_\alpha) = 1$ and is such that 
\[
h_{\rm{top}}(C^{(N)}_{\alpha}) = h_{\nu_{N}}(T). 
\]	
By assumption $\rm{(A3)}$, the function $\sum_{j\ge 1}q_{j}^{(N)}(\phi_j-\alpha_j)$ is well defined and has a unique equilibrium state. Now extending each $q^{(N)}$ putting zero in all coordinates $q_j^{(N)}$ with $j > N$, by compactness (Tychonoff's theorem) one can extract a subsequence of vectors $q^{(N_k)}$ in $Q_R$ converging to some $q^{*} \in Q_R$. 

Write
$\psi^{N}:=\sum_{j\le N}q^{(N)}_{j}(\phi_{j}-\alpha_{j})$. For any fixed
$M$ and any $N_{k}\ge M$, we get
\[
\begin{split}
\bigl\|\psi^{N_{k}}-\psi_{q^{*}}\bigr\|_{\infty}
&\le\Bigl\|\sum_{j\le M}\bigl(q^{(N_{k})}_{j}-q^{*}_{j}\bigr)(\phi_{j}-\alpha_{j})\Bigr\|_{\infty}
	\\
&+\Bigl\|\sum_{M<j\le N_{k}}q^{(N_{k})}_{j}(\phi_{j}-\alpha_{j})\Bigr\|_{\infty}
+\Bigl\|\sum_{j>N_k}q^{*}_{j}(\phi_{j}-\alpha_{j})\Bigr\|_{\infty} \\
& \le \Bigl\|\sum_{j\le M}\bigl(q^{(N_{k})}_{j}-q^{*}_{j}\bigr)(\phi_{j}-\alpha_{j})\Bigr\|_{\infty} + 2\varepsilon_{M}.
\end{split}
\]
On the right side of the above inequality, the first term is a finite sum whose coefficients tend to zero, so it vanishes as $k\to\infty$ with $M$ fixed. Hence 
\[
\limsup_{k \to \infty}\|\psi^{N_{k}}-\psi_{q^{*}}\|_{\infty}
\le2\varepsilon_{M} \quad \textrm{for every $M$}.
\]
Now letting $M\to\infty$, by assumption $\rm{(A2)}$ we finally get $\|\psi^{N_{k}}-\psi_{q^{*}}\|_{\infty}\to0$.

Invoking assumption $\rm{(A3)}$, let $\nu_{\alpha}$ be the unique equilibrium state for $\psi_{q^*}$, which is also ergodic. Since $\psi^{N_k}$ converges uniformly to $\psi_{q^*}$, and each $\psi^{N_k}$ has $\nu_{N_k}$ as its unique equilibrium state, it is well known that $\nu_{N_k}\to \nu_{\alpha}$ in the weak$^{*}-$ topology\footnote{In the literature, this property is mostly known as stability of equilibrium states.} (Theorem 4.2.11 in \cite{Kel98}). As a consequence, for each $j \ge 1$ we have
\[
\alpha_j = \lim_{k \to \infty}\int_{X}\phi_j\,d\nu_{N_k} = \int_{X}\phi_j\,d\nu_{\alpha},
\]
and hence $\nu_{\alpha} \in \cM(\alpha)$. Furthermore, by the ergodicity of $\nu_{\alpha}$, for each $j \ge 1$ the Birkhoff's ergodic theorem gives in particular
\[
\lim_{N \to \infty}\frac{1}{N}S_{N}\phi_j(x) = \int_{X}\phi_j\,d\nu_{\alpha} = \alpha_j
\]
for $\nu_{\alpha}-$almost every $x \in X$, that is, $\nu_{\alpha}(C_{\alpha_j}) = 1$ for each $j \ge 1$, where 
\[
C_{\alpha_j} := \biggl\{x \in X: \lim_{N \to \infty}\frac{1}{N}S_N\phi_j(x) = \alpha_j\biggr\}. 
\]
This readily implies that $\nu_{\alpha}$ is concentrated on the respective level set, that is, 
\[
\nu_{\alpha}(C_{\alpha}) = \nu_{\alpha}\bigg(\cap_{j = 1}^{\infty}C_{\alpha_j}\bigg) = 1. 
\]

Now we need to check what happens asymptotically with the metric entropy. Since each $\nu_{N_k}$ is the unique equilibrium state for $\psi^{N_k}$ and $\nu_{\alpha}$ is the unique equilibrium state for $\psi_{q^*}$, we can write
\[
\begin{split}
\big|h_{\nu_{\alpha}}(T) - h_{\nu_{N_k}}(T)\big| &\le \big|P_{\rm{top}}(\psi^{N_k}) - P_{\rm{top}}(\psi_{q^*})\big| + \bigg|\int_{X}\psi^{N_k}\,d\nu_{N_k} - \int_{X}\psi_{q^*}\,d\nu_{\alpha}\bigg|\\
& \le \big|P_{\rm{top}}(\psi^{N_k}) - P_{\rm{top}}(\psi_{q^*})\big| + \bigg|\int_{X}\psi^{N_k}\,d\nu_{N_k} - \int_{X}\psi_{q^*}\,d\nu_{N_k}\bigg| \\
& + \bigg|\int_{X}\psi_{q^*}\,d\nu_{N_k} - \int_{X}\psi_{q^{*}}\,d\nu_{\alpha}\bigg| \\
& \le 2\big\|\psi^{N_k}-\psi_{q^*}\big\|_{\infty} + \bigg|\int_{X}\psi_{q^*}\,d\nu_{N_k} - \int_{X}\psi_{q^{*}}\,d\nu_{\alpha}\bigg|,
\end{split}
\]
where we used that the topological pressure $\phi \mapsto P_{\rm{top}}(\phi)$ is a Lipschitz function (we assumed $h_{\textrm{top}}(T) < \infty$). Letting $k \to \infty$, we end up with $h_{\nu_{N_k}}(T) \to h_{\nu_{\alpha}}(T)$. Therefore, by item 1 in Theorem \ref{ACVP1}, we finally get 
\[
h_{\rm{top}}(C_{\alpha}) = \lim_{k \to \infty}h_{\nu_{N_k}}(T) = h_{\nu_{\alpha}}(T).
\]

We conclude the proof showing the uniqueness of the measure $\nu_{\alpha}$ supported on the level sets. Given $\alpha \in \cW(R)$ for some $R > 0$, suppose there exists $\eta_{\alpha} \in \cM(\alpha)$ with $\eta_{\alpha}(C_{\alpha}) = 1$ and such that $h_{\eta_{\alpha}}(T) = h_{\rm{top}}(C_{\alpha})$. In particular, we have
\[
\begin{split}
h_{\eta_{\alpha}}(T) + \int_{X}\psi_{q^*}\,d\eta_{\alpha} &= h_{\eta_{\alpha}}(T) + 0  \\
& = h_{\rm{top}}(C_{\alpha}) = h_{\nu_{\alpha}}(T) + 0 \\
& = h_{\nu_{\alpha}}(T) +  \int_{X}\psi_{q^*}\,d\nu_{\alpha} = P_{\rm{top}}(\psi_{q^*}),
\end{split}
\]
and then $\eta_{\alpha}$ is an equilibrium state for $\psi_{q^*}$. This implies $\eta_{\alpha} = \nu_{\alpha}$, and the theorem follows.
\end{proof}

\begin{remark}
Making the obvious changes, Theorems \ref{ACVP1} and \ref{ACVP2} (and other previous results) can be stated and proved in more generality, already considering the topological pressure instead of the entropy. However, to simplify the notation (as the entropy spectrum already captures the new phenomena with countable observables), we refrained to do so. In addition, for differentiable conformal systems and subshifts of finite type, Theorems \ref{ACVP1} and \ref{ACVP2} can be obtained for the Hausdorff dimension (for instance via the \emph{$u$-dimension} introduced in \cite{BS00}). 
\end{remark}
	
Especially because of Theorem \ref{ACVP2}, one might wonder now about the relationship between the interior points of rotation sets and the uniformly attainable points $\cW$. It is easy to see that in finite dimensions, when the observable $\Psi$ is $\R^d-$valued, under some suitable regularity assumptions, we have $\interior \rm{Rot}$$(\Psi) \subset \cW$, where the rotation set is calculated with respect to some types of permissible systems. On the other hand, the situation for $\R^{\N}-$valued observables is way more nuanced. To see that, we need some key results from the so-called \emph{Aubry-Mather theory} and \emph{ergodic optimization}. 
	
For a given continuous function $f\colon X \to \R$, write
\[
\beta_{\rm{max}}(f):= \max_{\mu\in\cM(T)}\int_X f\,d\mu, \quad \cM_{\rm{max}}(f) := \biggl\{\mu \in \cM(T): \int_{X}f\,d \mu = \beta_{\rm{max}}(f)\biggr\},
\] 
where the maximum is always attained because $\cM(T)$ is compact and $\mu \mapsto \int_{X}f\, d\mu$ is a continuous map. The \emph{Mather set} of $f$ is given by 
\[
\textrm{\rm{Ma}}(f):= \overline{\bigcup_{\mu \in \cM_{\rm{max}}(f)}\rm{supp}\,\mu},
\]
where $\rm{supp}\, \mu$ is the support of $\mu$. 
	
A continuous function $h\colon X \to \R$ is said to be a \emph{subaction} for $f$ if  
\[
f(x) + h(x) - (h\circ T)(x) \le \beta_{\rm{max}}(f) \quad \textrm{for all $x \in X$}. 
\]
A \emph{Mañé–Conze–Guivarc'h lemma} is a type of result that gives conditions on $T$ and $f$ for the existence of subactions. See for example \cite{Jen19} for a survey on ergodic optimization, and \cite{GL08} for the connections with the Aubry-Mather formulation.

\begin{proposition}\label{OPT}
Let $T$ be a subshift of finite type, an Anosov diffeomorphism, or a uniformly expanding map of a compact metric space $X$, which is also topologically mixing, and assume that each coordinate function $\psi_j$ is H\"older. Suppose that for every $N \ge 1$, except for the null function, there is no element in the space $\textrm{\rm{span}}\{\psi_1,\dots,\psi_N\}$ that is cohomologous to a constant. Then
\[
\cW \subset \interior_{\rm{ps}}\rm{Rot}(\Psi).
\]
\end{proposition}
\begin{proof}
Let $\alpha\in\cW$, that is, there exists $R > 0$ such that $\alpha \in \cW(R)$. Fix $N$ and $w = (w_1,\dots,w_N) \in\R^{N}\setminus\{0\}$, and consider the function $f:=\sum_{j\le N}w_{j}\psi_{j}$, which is also H\"older. Denote by $\nu:=\nu^{N}_{\alpha}$ the equilibrium state attaining $\alpha^{(N)}\in\cW_{N}(R)$, so that $\langle w,\alpha^{(N)}\rangle=\int_X f\,d\nu$. In this setup, it is well known that $\nu$ is a Gibbs measure (see the definition \eqref{GS}), consequently having a full support on $X$. 
		
Suppose, by contradiction, that 
\[
\int_X f\,d\nu=\beta_{\rm{max}}(f).
\]
Then $\nu \in \cM_{\rm{max}}(f)$, and since $\nu$ has full support, we immediately get $\rm{Ma}$$(f)=X$. Furthermore, by the Mañé–Conze–Guivarc'h lemma (\cite{Sav99,Bou01,CLT01,LT03}), $f$ admits a H\"older subaction, that is, there is a H\"older function $h$ with $f+h-h\circ T\le \beta_{\rm{max}}(f)$ everywhere. The equality is achieved on $\rm{Ma}$$(f)$. In fact, taking $g:= f + h - h\circ T - \beta_{\rm{max}}(f)$, for any $\mu \in \cM_{\rm{max}}(f)$ we get 
\[
\int_{X}g\,d\mu = \int_{X}f\,d\mu - \beta_{\rm{max}}(f) = 0.  
\]
But since $g \le 0$ on $X$, one can check that $g(x) = 0$ for every $x \in \textrm{supp}\,\mu$. By the arbitrariness of $\mu$, $g$ vanishes on the support of every maximizing measure for $f$, which, by continuity, also forces $g$ to vanish on $\rm{Ma}$$(f)$, as desired. 
		
Thus, the full support of $\nu$ implies $f$ is cohomologous to the constant $\beta_{\rm{max}}(f)$ everywhere, a contradiction. The same argument applied to $-f$ excludes the minimum. Hence, $\langle w,\alpha^{(N)}\rangle$ is strictly between the extreme values $\beta_{\rm{min}}(f)$ and $\beta_{\rm{max}}(f)$, where $\beta_{\rm{min}}(f) := \min_{\mu \in \cM(T)} \int_{X}f\,d\mu$. By the arbitrariness of $w$, this means
\begin{equation}\label{AWB}
\beta_{\rm{min}}(f) < \langle w,\alpha^{(N)} \rangle < \beta_{\rm{max}}(f) \quad \textrm{for all $w \in \R^{\N}\setminus \{0\}$}. 
\end{equation}
If $\alpha^{(N)} \in \partial \rm{Rot}$$(\psi^{(N)})$, the convexity of $\rm{Rot}$$(\psi^{(N)}) \subset \R^{N}$ guarantees the existence of $w' \neq 0$ such that $\langle w', y^{(N)}\rangle \le \langle w', \alpha^{(N)} \rangle$ for all $y^{(N)} \in \rm{Rot}$$(\psi^{(N)})$, which is a contradiction to \eqref{AWB}. Since $N$ is arbitrary, the result follows.   
\end{proof}
	
\begin{remark}
When $\Psi$ takes values in $\R^d$ ($d < \infty$), Proposition \ref{OPT} gives that $\cW = \interior \rm{Rot}(\Psi)$. This means that, for uniformly hyperbolic systems, the interior of the rotation set is the exact optimal place to look for ergodic full measures via regularity of the topological pressure.  
\end{remark}
	
\subsection{A geometric condition versus uniform attainability}\label{GCUA}
As observed above, for many permissible systems, interior points of rotation sets in $\R^d$ are always uniformly attainable via differentiability of the pressure function (see the condition \rm{(P)} below). To have, perhaps, a more direct way of checking if a parameter admits a full measure, we are tempted to extend a suitable idea of interior for the infinite-dimensional framework. 
	
Let $\alpha \in \rm{Rot}(\Psi)$, $N \ge 1$ and set 
\[
\delta_N(\alpha) := \textrm{dist}(\alpha^{(N)}, \partial \textrm{Rot}(\psi^{(N)})),
\]
where the distance is calculated with respect to some norm on $\R^N$. It is clear that $\delta_{M}(\alpha) \ge \delta_{N}(\alpha)$ whenever $M \le N$, and we then define the \emph{depth} of $\alpha$ in the rotation set as the number
\[
\delta(\alpha) = \lim_{N \to \infty}\delta_N(\alpha) = \inf_{N \ge 1}\delta_N(\alpha).
\]
We also consider the set of $\alpha$ with strictly positive depth, that is,
\[
\cA_{\rm{D}}:= \big\{\alpha \in \rm{Rot}(\Psi): \delta(\alpha) > 0 \big\}. 
\]
We call these points \emph{asymptotically deep}. The reason for this name is geometric. If $\delta(\alpha) > 0$, for each $N$ we have $\delta_N(\alpha) \ge \delta(\alpha) > 0$. This means that each $\rm{Rot}$$(\psi^N)$ contains a ball of radius $\delta_N(\alpha)$ and centered at $(\alpha_1,\dots,\alpha_N)$. When $N \to \infty$, $\delta_N(\alpha)$ may decrease but never reaches the border of $\rm{Rot}$$(\psi^{(N)})$.   
	
By the definition of depth and Proposition~\ref{OPT}, for many permissible systems and with $\Psi$ being $\R^d-$valued, it is clear that
\[
\cA_{\rm{D}} = \interior \textrm{Rot}(\Psi) = \cW. 
\]
In this case, we always have
\[
\textrm{interior vectors} \implies \textrm{differentiability of the pressure} \implies \textrm{uniform attainability}.
\]
By what we shall see in the rest of this work, the scenario is quite different for $\R^{\N}-$valued observables: in the above equivalences we have proper inclusions. More especifically, for many permissible systems, we have
\[
\cA_{\rm{D}} \subset \cW \subset \interior_{\rm{ps}}\rm{Rot}(\Psi).
\]
As mentioned earlier, for each $N \ge 1$ we work under the following standing assumption on the $\R^{N}-$valued observable $\psi^{(N)}$:
\begin{itemize}
\item[(P)] For every $N\ge1$ and $q\in\R^{N}$ the function $\sum_{j\le N}q_j\psi_j$ has a unique equilibrium state $\nu_q^{N}$, and the pressure function $P_N(q):=P_{\rm{top}}\bigl(\sum_{j\le N}q_j\psi_j\bigr)$ is differentiable on $\R^{N}$ with 
\[
\nabla P_N(q)=\biggl(\int_X\psi_1\,d\nu^{N}_q,\dots,\int_X\psi_N\,d\nu^{N}_q\biggr).
\]
\end{itemize}
All the permissible systems described in Remark \ref{PSYST}, together with any observable with coordinate functions belonging to any vector space of functions admitting unique equilibrium states (the set $\cU$ in the definition of permissible systems is the natural candidate), satisfy the property \rm{(P)}. For these systems, it makes sense to compare asymptotic depth with uniform attainability. 
\begin{proposition}\label{PP}
Assume $\rm{(P)}$. Then
\[
\alpha \in \cA_{\rm{D}} \implies \alpha \in \cW\biggl(\frac{h_{\rm{top}}(T)}{\delta(\alpha)}\biggr).
\]
In particular, $\cA_{\rm{D}} \subseteq \cW$. 
\end{proposition}
\begin{proof}
Fix $N$ and consider the convex function $G_N(q):=P_N(q)-\langle q,\alpha^{(N)}\rangle$ on $\R^{N}$. Given $q\neq0$, the point $y :=\alpha^{(N)}+\delta_N(\alpha)\,q/\|q\|_2$ lies in $\rm{Rot}$$(\psi^{(N)})$ (recall that $\|\cdot\|_{\infty} \le \|\cdot\|_{2}$). Hence $y =\int\psi^{(N)}d\mu$ for some $\mu\in \cM(T)$, and the variational principle gives
\[
G_N(q)\ \ge\ h_\mu(T)+\langle q,y\rangle-\langle q,\alpha^{(N)}\rangle
\ \ge\ \delta_N(\alpha)\,\|q\|_2\ \ge\ \delta(\alpha)\,\|q\|_2 .
\]
This implies that $G_N(q)$ takes arbitrarily large values when $\|q\|_{2} \to \infty$. In particular, there exists $L > 0$ such that $\|q\|_{2} \ge L$ gives $G_N(q) \ge G_N(0)$. 
Then, the minimum of $G_N$ is attained at some point $q^{*}$, and we have
\[
\delta(\alpha)\,\|q^{*}\|_2\ \le\ G_N(q^{*})\ \le\ G_N(0)=P_N(0)= h_{\rm{top}}(T),
\]
which immediately yields $\|q^{*}\|_\infty\le\|q^{*}\|_2\le h_{\rm{top}}(T)/\delta(\alpha)$. Moreover, by the property $\textrm{(P)}$ the function $G_N$ is differentiable and must satisfy 
\[
0=\nabla G_N(q^{*})= \biggl(\int\psi_1\,d\nu^{N}_{q^{*}},\dots,\int\psi_N\,d\nu^{N}_{q^{*}}\biggr)-\alpha^{(N)},
\]
where $\nu^{N}_{q^{*}}$ is the unique equilibrium state of $\sum_{j\le N}q^*_j\psi_j$. Therefore, for every $N \ge 1$ we get $\alpha^{(N)} \in \cW_{N}(h_{\rm{top}}(T)/\delta(\alpha))$, and the result follows. 
\end{proof}
	
In the following, we see that positive asymptotic depth forces positive \emph{oscillation} in each coordinate function of the observable, which is a fatal obstruction to producing ergodic full measures via the mechanism of Theorem~\ref{ACVP2}.  
	
\begin{proposition}\label{OSC}
If $\delta(\alpha) > 0$ then
\[
\textrm{$\rm{osc}$}(\psi_j) := \max_{x \in X}\psi_j(x) - \min_{x \in X}\psi_j(x) \ge 2\delta(\alpha) \quad \textrm{for every $j \ge 1$}. 
\]
In particular, \textup{(A2)} fails on $\cA_{\rm{D}}$.  
\end{proposition}
\begin{proof}
Fix $j \ge 1$ and take $N\ge j$. The points $\alpha^{(N)}\pm\delta(\alpha)e_j$ are in $\overline{B(\alpha^{(N)},\delta_N(\alpha))}\subseteq\rm{Rot}$$(\psi^{(N)})$, and then there are $\mu_{\pm}\in \cM(T)$ with $\int\psi_j\,d\mu_{\pm}=\alpha_j\pm\delta(\alpha)$.
Hence the range of $\mu\mapsto\int\psi_j\,d\mu$ has length at least $2\delta(\alpha)$, and this range is contained in $[\min\psi_j,\max\psi_j]$, giving $\rm{osc}(\psi_j)\ge2\delta(\alpha)$ and $2\|\psi_j\|_\infty\ge \rm{osc}(\psi_j)$. 
		
Let us now show how the series diverge. For each fixed $R>0$ and $N \ge 1$, set
\[
\varepsilon_N:=\sup_{\sup_j|\gamma_j|\le R}\ \Bigl\|\sum_{j>N}
\gamma_j(\psi_j-\alpha_j)\Bigr\|_\infty,
\]
and take $\gamma = R e_m$, with $m > N$. In this case, we get
\[
\varepsilon_N \ge R\|\psi_m-\alpha_m\|_\infty\ge \frac{1}{2}R\,\textrm{osc}(\psi_m)\ge R\,\delta(\alpha).
\]
Thus, we do not have $\varepsilon_N \to 0$, which implies the divergence of $\sum_{j \ge 1}\|\psi_j\|_{\infty}$. 
\end{proof}
	
The picture is now quite clear: to admit an ergodic full measure via the mechanism of Theorem \ref{ACVP2}, $\alpha \in \rm{Rot}(\Phi)$ must belong to $\cW \setminus \cA_{\rm{D}}$. On the other hand, if assumption \textup{(A2)} of Theorem~\ref{ACVP2} holds for some $\alpha\in\rm{Rot}(\Phi)$ and some $R>0$, then $\textrm{$\rm{osc}$}(\phi_j)\to0$ and $\cA_{\rm D}=\emptyset$. In other words, for a given observable $\Phi$, either $\cA_{\rm D}=\emptyset$ or Theorem~\ref{ACVP2} applies to no vector in $\rm{Rot}(\Phi)$. 
    
Interestingly, we did not find any pair $(T,\Phi)$ satisfying $\cA_{\rm{D}} \neq \emptyset$ in the rest of this work. The construction of such examples seems to be very convoluted. Nevertheless, understanding $\cA_{\rm{D}}$ would bring a more geometric flavor to the infinite-dimensional rotation sets, in the sense of Kucherenko-Wolf results \cite{KW14}. 
 	
\subsection{No ergodic full measure for a parameter outside $\cW$}\label{NW}
	
Let $\Sigma = \{0,1\}$, consider the shift map $\sigma\colon \Sigma^{\N} \to \Sigma^{\N}$. Considering $s:= s(\omega,\widetilde{\omega})$ on $\Sigma^{\N} \times \Sigma^{\N}$ as the smallest positive number $s$ such that $\omega_{s} \neq \widetilde{\omega}_{s}$, we define the distance on $\Sigma^{\N}$ to be $d(\omega, \widetilde{\omega}) = 2^{-s(\omega,\widetilde{\omega})}$ if $\omega \neq \widetilde{\omega}$ and $d(\omega, \widetilde{\omega}) = 0$ if $\omega = \widetilde{\omega}$. Let $\Psi\colon \Sigma^{\N} \to \R^{d}$ with all coordinate functions $\psi_j$ being Lipschitz on $\Sigma^{\N}$. In this finite-dimensional situation, it is well known that each $\alpha \in \interior \rm{Rot}(\Psi)$ admits an ergodic full measure (\cite{BSS02b}). In the following construction, we are going to show that this might not be the case when $\Psi$ has a countable number of locally constant (in particular Lipschitz) coordinate functions but $\alpha$ is not uniformly attainable.
	
\begin{theorem}\label{NEFM}
Given $\epsilon > 0$ there exists an observable $\Phi\colon \Sigma^{\N} \to \R^{\N}$ with locally constant coordinate functions $\phi_j$ and $\alpha \in \rm{Rot}(\Phi)$ such that:
\begin{enumerate}
\item $\sum_{j \ge 1}\|\phi_j\|_{\infty} < \infty$ and $\sum_{j \ge 1}\rm{Lip}$$(\phi_j) < \infty$.
			
\item $\alpha \in \interior_{\rm{ps}}\rm{Rot}(\Phi)$.
			
\item $\alpha \notin \cW$. 
			
\item $h_{\rm{top}}(C_{\alpha}) > h_{\rm{top}}(\sigma) - \epsilon$.   
			
\item $\mu(C_\alpha) = 0$ for every $\sigma-$invariant measure $\mu$. 
\end{enumerate}
In particular $\alpha$ does not admit any full measure, ergodic or not.  
\end{theorem}
	
The proof of this result is, in part, inspired by \cite{Bou01} and will be divided into some lemmas. 
	
Let $C^0(\Sigma^{\N},\R)$ be the space of real-valued continuous functions on $\Sigma^{\N}$, and $L_n\subset C^0(\Sigma^{\N},\R)$ the space of functions depending only on the first $n$ coordinates, so $\dim L_n=2^n$. $L_0$ is the space of constants, and $L:=\bigcup_{n\ge0}L_n$ is uniformly dense in $C(\Sigma^{\N},\R)$. We set
\[
E_0:=\R,\qquad
E_n:=\R +\{g\circ\sigma-g:\ g\in L_{n-1}\}\subseteq L_n \quad \textrm{for $n \ge 1$}. 
\]
Note that $g\in L_{n-1}$ implies $g\circ\sigma\in L_n$, so indeed $E_n\subseteq L_n$.
	
\begin{lemma}\label{lem:algebra}
For every $n\ge1$, we have 
\begin{enumerate}[label=\textup{(\alph*)},leftmargin=2.6em]
\item $\dim E_n=2^{n-1}$.
\item $E_{n-1}\subseteq E_n$ and $E_n\cap L_{n-1}=E_{n-1}$.
\end{enumerate}
\end{lemma}
	
\begin{proof}
For the first item, consider $D\colon L_{n-1}\to L_n$, $D(g):=g\circ\sigma-g$. If $D(g)=0$ then $g$ is a continuous $\sigma$-invariant function, hence constant on the dense orbit of a transitive point, and thus $g$ is forced to be constant everywhere. Then, we get $\ker D=\R$ and $\dim D(L_{n-1})=2^{n-1}-1$. Moreover, $\R\cap D(L_{n-1})=\{0\}$. In fact, if $g\circ\sigma-g=c$, integrating with respect to any $\mu\in\M(\sigma)$ gives $c=0$. Hence, $\dim E_n=1+(2^{n-1}-1)=2^{n-1}$.
		
Now let us prove item (b). For $n\ge1$ we have $L_{n-2}\subseteq L_{n-1}$ (with $L_{-1}:=L_0$), which comes from $E_{n-1}\subseteq E_n$. We also have $E_{n-1}\subseteq L_{n-1}$, and this gives $E_{n-1}\subseteq E_n\cap L_{n-1}$. For the reverse inclusion, the case $n=1$ is immediate: $L_0$ consists of constants, so $E_1=\R=E_0$ and $E_1\cap L_0=E_0$. Let $n\ge2$, let $g\in L_{n-1}$ and $c\in\R$, and suppose $f:=c+g\circ\sigma-g$ lies in $L_{n-1}$. Writing $g=g(x_1,\dots,x_{n-1})$ we have $(g\circ\sigma)(x)=g(x_2,\dots,x_n)$. Both $f$ and $g$ are independent of $x_n$, hence so is $g\circ\sigma=f-c+g$; that is, $g(x_2,\dots,x_n)$ does not depend on $x_n$, which means that $g$ does not depend on its last argument and $g\in L_{n-2}$. Therefore $f\in E_{n-1}$.
\end{proof}
	
\begin{lemma}\label{lem:complement}
There are finite sets $\cB_n\subset L_n$, with $\#\cB_1=1$ and $\#\cB_n=2^{n-2}$ for $n\ge2$ such that, setting $V_0:=\{0\}$ and $V_n:=\Span(\cB_1\cup\dots\cup\cB_n)$,
\[
L_n=E_n\oplus V_n \qquad\text{for every } n\ge1 .
\]
In particular $\dim V_n=2^{n-1}$, $V_{n-1}\subseteq V_n$, and $E_m\cap V_n=\{0\}$ whenever $m\ge n\ge1$. Moreover $\cB_1\cup\cB_2\cup\cdots$ is a linearly independent set.
\end{lemma}
	
\begin{proof}
We use induction on $n$. For $n=1$, $\dim L_1=2$ and $E_1=\R$, so any $\cB_1=\{\varphi\}$ with $\varphi\in L_1\setminus\R$ works. Assume $L_{n-1}=E_{n-1}\oplus V_{n-1}$. Then
\[
E_n\cap V_{n-1}\ \subseteq\ E_n\cap L_{n-1}\cap V_{n-1}
\ =\ E_{n-1}\cap V_{n-1}\ =\ \{0\}
\]
by Lemma~\ref{lem:algebra}(b), so the sum $E_n+V_{n-1}$ is direct and has dimension $2^{n-1}+2^{n-2}$. Choose $\cB_n\subset L_n$ consisting of $2^n-2^{n-1}-2^{n-2}=2^{n-2}$ vectors completing $E_n\oplus V_{n-1}$ to $L_n$. Then, $L_n=E_n\oplus V_n$ with $V_n=V_{n-1}\oplus\Span\cB_n$ and $\dim V_n=2^{n-1}$. Finally, for $m\ge n$ we have $V_n\subseteq V_m$ and $E_n\subseteq E_m$, so $E_m\cap V_n\subseteq E_m\cap V_m=\{0\}$. Linear independence of $\bigcup_n\cB_n$ is immediate from $V_{n}=V_{n-1}\oplus\Span\cB_n$.
\end{proof}
	
Enumerate the set $\bigcup_{n\ge1}\cB_n=\{\phi^0_1,\phi^0_2,\dots\}$ in any order, set
\begin{equation}\label{re-scale}
\phi_j:=c_j\,\phi^0_j \quad \textrm{with} \quad 
c_j:=\frac{2^{-j}}{1+\|\phi^0_j\|_\infty+\rm{Lip}(\phi^0_j)},
\end{equation}
and consider the observable $\Phi$ with each coordinate given by $\phi_j$. Note that every $f \in L_n$ is Lipschitz. In fact, if $d(x,y)\le2^{-(n+1)}$ then $f(x)=f(y)$, and otherwise $|f(x)-f(y)|\le2\|f\|_\infty\le 2^{n+1}\|f\|_\infty\,d(x,y)$. Thus $\rm{Lip}$$(f)\le2^{n+1}\|f\|_\infty<\infty$.
	
Take any $p \in (0,1)$ with $p \neq 1/2$, and finally define
\[
\eta:=\tfrac12\bigl(m_{p}+m_{1-p}\bigr),
\qquad \alpha_j:=\int_{\Sigma^{\N}}\phi_j\,d\eta,
\]
where $m_p$ is the \emph{Bernoulli measure} with $m_p[1]=p$, and the set $[1]$ is taken as the cylinder formed of points $(w_0 w_1 \dots)$ with $w_0 = 1$. In particular, note that $m_p \neq m_{1-p}$, meaning that $\eta$ is not an extreme point of $\M(\sigma)$, and hence is not ergodic. 
        
\begin{lemma}\label{lem:sep}
The integration map $\M(\sigma)\ni\mu\mapsto\bigl(\int_{\Sigma^{\N}}\phi_1\,d\mu, \int_{\Sigma^{\N}}\phi_2\,d\mu,\dots\bigr) \in \R^{\N}$ is injective, that is, it separates measures. In particular, $\M(\alpha)=\{\eta\}$. 
\end{lemma}
	
\begin{proof}
Let $\widehat \alpha$ be arbitrary, take $\mu,\nu\in\M(\widehat \alpha)$, and let $f\in L$, say $f\in L_n$. It follows from Lemma~\ref{lem:complement} and definition \eqref{re-scale} that 
\[
f=c + g\circ\sigma-g +\sum_{j\in F}b_j\phi_j
\]
with $c\in\R$, $g\in L_{n-1}$, $F$ finite and each $b_j\in\R$. Since the measures are invariant, we get 
\[
\int_{\Sigma^{\N}}f\,d\mu=c+\sum_{j\in F}b_j\int_{\Sigma^{\N}}\phi_j\,d\mu \quad \textrm{and} \quad \int_{\Sigma^{\N}} f\,d\nu=c+\sum_{j\in F}b_j\int_{\Sigma^{\N}}\phi_j\,d\nu,
\]
which readily implies $\int_{\Sigma^{\N}}f\,d\mu = \int_{\Sigma^{\N}} f\,d\nu$ for every $f\in L$. As the space $L$ is uniformly dense in $C^{0}(\Sigma^{\N},\R)$, we must have $\mu=\nu$. The arbitrariness of $\widehat{\alpha}$ guarantees the injectivity of the map, and of course we obtain $\cM(\alpha) = \{\eta\}$. 
\end{proof}
	
\emph{Proof of Theorem \ref{NEFM}}. By the definition of the observable in \eqref{re-scale}, we have
\[
\|\phi_j\|_{\infty} \le \frac{2^{-j}\|\phi^{0}_j\|_{\infty}}{1 + \|\phi^{0}_j\|_{\infty} + \textrm{Lip}(\phi^{0}_j)} \le 2^{-j} \quad \textrm{for all $j \ge 1$.}
\]
Similarly, we have the bound
\[
\textrm{Lip}(\phi_j) \le c_j\textrm{Lip}(\phi^{0}_j) = \frac{2^{-j}\textrm{Lip}(\phi^{0}_j)}{1 + \|\phi^{0}_j\|_{\infty} + \textrm{Lip}(\phi^{0}_j)} \le 2^{-j} \quad \textrm{for all $j \ge 1$.}
\]
In particular, $\Phi$ satisfies item (1).  
	
Let us obtain item (2). Fix $N \ge 1$, take $b = (b_1,\dots, b_N) \in \R^{N}\setminus \{0\}$, and consider the function $f := \sum_{j\le N}b_j \phi_j$. Note that $f\ne0$ by the linear independence in Lemma~\ref{lem:complement}. Moreover, it is clear that $f\in V_n\subseteq L_n$ for some $n \ge 1$. Suppose, by contradiction, that 
\[
\int_{\Sigma^{\N}}f\,d\eta= \max_{\mu\in\cM(\sigma)}\int_{\Sigma^{\N}}f\,d\mu =: \beta_{\rm{max}}(f). 
\]
By the Mañé–Conze–Guivarc'h lemma for locally constant functions, there exists $h\in L_{n-1}$ such that 
\[
f + h - h\circ\sigma \le \beta_{\rm{max}}(f).  
\]
Letting $u:= h\circ\sigma-h - f + \beta_{\rm{max}}(f) \ge\ 0$, and integrating with respect to $\eta$, we obtain 
\[
\int_{\Sigma^{\N}}u\,d\eta = -\int_{\Sigma^{\N}}f\,d\eta + \beta_{\rm{max}}(f) = 0. 
\]
Since $u \ge 0$, one can check that $u = 0$ on $\textrm{supp}\, \eta$. Since $\eta$ has full support and $u$ is continuous, $u$ is forced to be zero everywhere. This gives
\[
f= h\circ\sigma-h + \beta_{\rm{max}}(f) \in E_n. 
\]
But $f\in V_n$ and $E_n\cap V_n=\{0\}$ by Lemma~\ref{lem:complement}. Therefore, $f=0$ everywhere, a contradiction. Following in the same way for $-f$, we conclude that $\langle b,\alpha^{(N)}\rangle$ is strictly between the extreme values $\beta_{\rm{min}}(f)$ and $\beta_{\rm{max}}(f)$, with $\beta_{\rm{min}}(f) := \min_{\mu \in \cM(\sigma)} \int_{\Sigma^{\N}}f\,d\mu$. The arbitrariness of $b$ and convexity of $\rm{Rot}$$(\phi^{(N)})$ guarantee that for each $N \ge 1$, we have $\alpha^{(N)} \in \interior \rm{Rot}$$(\phi^{(N)})$, as desired.  
	
For item (3), assume that $\alpha \in \cW(R)$ for some $R > 0$. Invoking the stability of equilibrium states, exactly as in the proof of Theorem \ref{ACVP2}, $\eta$ must be an equilibrium state for a Lipschitz function, that is, $\eta$ is ergodic, a contradiction. Thus, $\alpha \notin \cW$.
	
To show item (4), we first claim that
\begin{equation}\label{eq:VC} C_\alpha=\{x\in\Sigma^{\N}:\ V(x)\subseteq\M(\alpha)\},
\end{equation}
where $V(x)$ is the set of limit points of $\nu_N := \frac{1}{N}\sum_{k=0}^{N-1}\delta_{\sigma^{k}(x)}$. In fact, for a fixed $j \ge 1$, the bounded sequence
\begin{equation}\label{SNP}
\frac{1}{N}S_N\phi_j(x) = \int_{\Sigma^{\N}}\phi_j \, d\nu_N(x) 
\end{equation}
converges to $\alpha_j$ if and only if all of its limit points are equal to $\alpha_j$. Since its set of limit points is exactly $\{\int_{\Sigma^{\N}}\phi_j\,d\mu:\mu\in V(x)\}$, $V(x) \subseteq \cM(\alpha)$ readily implies $x \in C_{\alpha}$. Conversely, let $x \in C_{\alpha}$. The equality in \eqref{SNP}, together with the compactness of $\M(\sigma)$, guarantees the existence of a subsequence $\nu_{n_k}(x)$ converging to some $\mu \in \cM(\sigma)$ such that 
\[
\alpha_j = \lim_{k \to \infty}\int\phi_j(x)\,d\nu_{n_k}(x) = \int\phi_j\,d\mu \quad \textrm{for each $j \ge 1$},
\]
that is, $V(x) \subset \cM(\alpha)$, and the claim follows. By Lemma~\ref{lem:sep}, $\M(\alpha)=\{\eta\}$, and then \eqref{eq:VC} becomes $C_\alpha=\{x:V(x)=\{\eta\}\}=G_\eta$. The saturation property now gives $h_{\rm{top}}(C_\alpha)=h_{\rm{top}}(G_\eta) =h_{\eta}(\sigma)$. In particular, recalling $\eta = (m_{p} + m_{1-p})/2$, and using the affinity of the metric entropy, we obtain
\[
h_{\eta}(\sigma) = \frac{1}{2}(H(p) + H(1-p)) = H(p) := -p\log p - (1-p)\log(1-p).
\] 
Thus, letting $p$ arbitrarily close to 1/2, implies $h_{\rm{top}}(C_{\alpha})$ arbitrarily close to $\log 2 = h_{\rm{top}}(\sigma)$, as desired.   
	
Finally, as the level set $C_{\alpha}$ is exactly the set of $\eta-$generic points, we immediately get $\mu(C_{\alpha}) = 0$ for every $\mu \in \cM(\sigma)\setminus \{\eta\}$. In particular, this implies item (5), and the theorem is proved. \qed  

\begin{remark}
For an explicit example of the observable $\Phi$ in Theorem \ref{NEFM}, consider $\mathcal U$ as the set of finite words beginning and ending with $1$, including the word $1$. Taking any enumeration $\{u_1, u_2,\dots\}$ of $\cU$, one can take $\Phi$ as the observable with coordinate functions given by 
\[
\phi_j = 8^{-|u_j|}\mathbf 1_{[u_j]}.
\] 
This example is a subfamily of the observable in Theorem~\ref{REG}. The words beginning or ending with $0$, which are discarded here, are exactly those responsible for the linear relations $\mathbf 1_{[w0]}+\mathbf 1_{[w1]}=\mathbf 1_{[w]}$ that make the pseudo-interior of the rotation set of that observable empty.
\end{remark}

\subsection{Ergodic full measure for a parameter outside $\cW$}\label{OTW}
	
Since the sets $\cA_{\rm{D}}$ and $\cW$ are equivalent to $\interior \rm{Rot}(\Phi)$ when $\Phi$ takes values in $\R^d$, to produce ergodic full measures outside $\cW$, we only need to check what happens to the points in $\partial \rm{Rot}(\Phi)$. However, in the infinite-dimensional framework, there is way more room to explore in $\rm{Rot}(\Phi)\setminus \cW$, but the ergodic full measures cannot be produced via Theorem \ref{ACVP2}. In this subsection, we show how to construct them through a different mechanism. 

Let $\Sigma = \{0,1\}$, and write $[\omega]:= \{x \in \Sigma^{\N}: (x_1 x_2 \dots x_{|\omega|}) = w\}$, the cylinder generated by a finite word $\omega$. For each $j \ge 1$, define
\[
A_j:= [0^{j-1} 1] \quad \textrm{and} \quad \xi_{j}:= \textbf{1}_{A_j}. 
\]
So, our observable is $\Xi: \Sigma^{\N} \to \R^{\N}$ given by $\Xi(x) = (\xi_1(x),\xi_2(x),\dots)$. Observe that $\cM(\alpha) = \{\mu \in \cM(\sigma): \mu(A_j) = \alpha_j \quad \textrm{for all $j \ge 1$}\}$. 

We now fix $\alpha = (\alpha_1,\alpha_2,\dots)$ as 
\[
\alpha_j = c4^{-2^j} \quad \textrm{with} \quad c:= \biggl(\sum_{j\ge 1}4^{-2^{j}}\biggr)^{-1}. 
\]
It is clear that $\alpha_j \in (0,1)$ for every $j \ge 1$, and $\sum_{j\ge 1}\alpha_j = 1$. \\ 

\begin{theorem}\label{NTW}
Let the observable $\Xi$ and the vector $\alpha$ be as above. Then 
\begin{enumerate}
\item $\alpha \in \interior_{\rm{ps}}\rm{Rot}(\Xi) \setminus \cW$. 

\item $\alpha$ admits an ergodic full measure, and $ 0 < h_{\rm{top}}(C_{\alpha}) < \log 2$. 
\end{enumerate}
\end{theorem} 

The proof is divided into several auxiliary steps. The main idea is to invoke an induced system where the induced measure is actually linked to a given a probability vector on $\N$, which also depends on the vector $\alpha$. Even though we are working on compact spaces, the construction uses some tools from Sarig's formalism for countable shifts (\cite{Sar99, Sar01}), a setup with a non-compact phase space. 

We start by introducing some needed machinery. Let $B:= A_1 = [1]$. For $x \in B$, define the \emph{return function} $r$ and the \emph{induced transformation} $\sigma_B$ as
\[
r(x) = \min\{n \ge 1: x_{n} = 1\} \quad \textrm{and} \quad \sigma_B(x):= \sigma^{r(x)}(x).
\]
We now take the set $B':= \{x \in B: x_{n} = 1 \textrm{ for infinitely many $n$}\}$. For each $x \in B'$, we define 
\[
\rho(x) = \bigl(r(x),r(\sigma_{B}(x)), r(\sigma_{B}^{2}(x)),\dots\bigr).  
\]
That is, a point in $B'$ is completely encoded by the sequence of distances between consecutive symbols $1$:
\[
B' \ni x = (1\,\underbrace{0\cdots0}_{k_1-1}1\,\underbrace{0\cdots0}_{k_2-1}1\cdots) \textrm{ is associated to } \rho(x)=(k_1,k_2,k_3,\dots) \in \N^{\N}.
\]
It is easy to see that $\rho: B' \to \N^{\N}$ is a bijection. Given a probability vector $(P_k)_{k \ge 1}$, let $\mu:=(\rho^{-1})_{*}\bigl(P^{\otimes\N}\bigr)$ be the push-forward of the Bernoulli measure $P^{\otimes\N}$ under the homeomorphism
$\rho$. This is a Borel probability on $\Sigma^{\N}$ concentrated on $B'$. More specifically, $\mu$ is the unique Borel probability with $\mu(B')=1$ and
\begin{equation}\label{CBS}
\mu\Bigl(\bigl[\,1\,0^{k_1-1}1\cdots1\,0^{k_n-1}1\,\bigr]\Bigr)
=\prod_{i=1}^{n} P_{k_i}.
\end{equation}

For any arbitrary $\ell \in \N$, we denote by $\{r = \ell\}$ the set of points $x \in \Sigma^{\N}$ such that $r(x) = \ell$. The following result is based on a classical construction known as \emph{Kakutani towers} \cite{Kak43, VO16}. It is the main auxiliary lemma in our approach. 

\begin{lemma}\label{RNW} Let $(P_k)_{k \ge 1}$ be a probability vector on $\N$ with $p:= \sum_{k \ge 1}kP_k < \infty$, and let $\overline{\nu}$ be the push-forward of the Bernoulli measure under the map $\rho^{-1}$. Then, the measure on $\Sigma^{\N}$  
\[
\nu(E) = \frac1{p}\sum_{\ell \ge 1}\sum_{k \le \ell -1}\overline{\nu}(\{r = \ell\} \cap \sigma^{-k}(E)) 
\]
is $\sigma-$invariant, ergodic and satisfies
\[
\nu(B) = \frac{1}{p}, \quad \nu(A_j) = \frac{1}{p} \sum_{\ell \ge j}P_{\ell} \quad \textrm{for every $j \ge 1$}.
\]
Moreover, the metric entropy of $\nu$ is given by 
\[
h_{\nu}(\sigma) = -\frac{1}{p}\sum_{k \ge 1}P_k \log P_k.  
\]
\end{lemma}

\begin{proof}
It is well known that $\nu$ is an ergodic $\sigma-$invariant probability measure (see for example Section 1.5 in \cite{CFS82}, Theorem 1). 

Fix $\ell\ge1$ and set
\[
j(k):=\begin{cases}1,&k=0,\\ \ell+1-k,&1\le k\le\ell-1.\end{cases}
\]
We claim that, for all $0\le k\le\ell-1$ and all $j\ge1$,
\begin{equation}\label{AON}
\{r=\ell\}\cap\sigma^{-k}A_j=
\begin{cases}
\{r=\ell\}, & j=j(k),\\[2pt]
\emptyset, & j\neq j(k),
\end{cases}
\end{equation}
and that $k\mapsto j(k)$ is a bijection from $\{0,1,\dots,\ell-1\}$ onto
$\{1,2,\dots,\ell\}$.

In fact, on $\{r=\ell\}$ the coordinates $x_0,\dots,x_{\ell}$ vanish except for $x_0=x_{\ell}=1$. Take first $k=0$. For $j\ge2$ the set $\sigma^{-k}(A_j)$ requires $x_0=0$, which is incompatible with $x_0=1$. For $j=1$ it only requires $x_0=1$, which is, of course, already true. Suppose now $1\le k\le\ell-1$, and distinguish three cases according to the position
$k+j-1$ of the symbol $1$ prescribed by $\sigma^{-k}(A_j)$. If $k+j-1\le\ell-1$,
that is $j\le\ell-k$, then $x_{k+j-1}=0$ on $\{r=\ell\}$ and the intersection
is empty. If $k+j-1=\ell$, that is $j=\ell+1-k$, then the conditions imposed by $\sigma^{-k}(A_j)$ force $x_k=\cdots=x_{\ell-1}=0$ and $x_{\ell}=1$, all of
which hold on $\{r=\ell\}$, meaning that the intersection is the whole set. If $k+j-1\ge\ell+1$, that is $j\ge\ell+2-k$, then $k+j-2\ge\ell$ and $\sigma^{-k}(A_j)$ requires $x_{\ell}=0$, contradicting $x_{\ell}=1$. This proves \eqref{AON}. In addition, $j(0)=1$, while $j(k)$ decreases from $\ell$ to $2$ as $k$ increases from $1$ to $\ell-1$. That is, the correspondence $k \mapsto j(k)$ is a bijection onto $\{1,\dots,\ell\}$ as desired. 

Consequently, each term $\bar\nu(\{r=\ell\}\cap\sigma^{-k}(A_j))$ equals $P_{\ell}$ or $0$, with precisely one $k$ contributing when $j\le\ell$ and none when $j>\ell$. Then, by the definition of $\nu$, for each $j \ge 1$ we finally obtain
\[
\nu(A_j)
=\frac{1}{p}\sum_{\ell\ge1}\sum_{k=0}^{\ell-1}\bar\nu\bigl(\{r=\ell\}\cap\sigma^{-k}(A_j)\bigr) = \frac{1}{p}\sum_{\ell\ge j}P_{\ell}.
\]

We now show the entropy formula. When $p = 1$, the result holds trivially. Suppose $p>1$ and let
$G_k:=(1-\theta)\theta^{k-1}$ with $\theta:=1-1/p$. It follows from Gibbs' inequality that
\[
-\sum_{k \ge 1} P_k\log P_k \le-\sum_{k \ge 1}P_k\log G_k=-\log(1-\theta)-(p-1)\log\theta<\infty .
\]
Now define the partition
\begin{equation}\label{PPR}
\mathcal{P}_r := \big\{\rho^{-1}[k]: k \in \N\big\},\quad \textrm{where} \quad [k]:=\{\omega\in\mathbb{N}^{\mathbb{N}}:\omega_0=k\}.
\end{equation}
Since $\{[k]: k \in \N\}$ is a generating partition of $\N^{\N}$ with entropy $-\sum_kP_k\log P_k < \infty$, and the transformation induced by $\sigma$ on $B$ is $\sigma_B$ with invariant measure $\bar\nu$, one can check that $\mathcal{P}_{r}$ is a generating partition of $B$. In this case, we have 
\[
H_{\bar\nu}(\mathcal{P}_r) = -\sum_{k \ge 1} \bar{\nu}(\{r= k\})\log \bar{\nu}(\{r= k\})<\infty,
\]
and, by Sinai's generator theorem (see Theorem 4.1 in \cite{Sar09}), we have $h_{\bar\nu}(\sigma_B) = h_{\bar{\nu}}(\sigma_B,\mathcal{P}_r) \le H_{\bar\nu}(\mathcal{P}_r)$. By using formula \eqref{CBS}, one can calculate $h_{\bar{\nu}}(\sigma_B,\mathcal{P}_r) = -\sum_{k \ge 1} P_k \log {P_k}$. Therefore, Abramov's formula (\cite{Abr59}) together
with $\nu(B)=1/p$ yields
\[
h_{\nu}(\sigma)=\nu(B)\,h_{\bar\nu}(\sigma_B) = \frac{1}{p}{\sum_{k\ge1} P_k\log P_k}. \qedhere
\]
\end{proof}

\begin{lemma}\label{ENT}
The measure $\nu_{\alpha}$ produced by Lemma \ref{RNW} with $P_{\ell} = (\alpha_{\ell} - \alpha_{\ell + 1})/\alpha_1$ is ergodic, satisfies $\nu_{\alpha}(A_j) = \alpha_j$ for every $j\ge1$ and $\nu_{\alpha}(C_{\alpha}) = 1$. Moreover, we have
\[
h_{\nu_{\alpha}}(\sigma) = \alpha_1\log \alpha_1 - \sum_{j \ge 1}(\alpha_j - \alpha_{j+1})\log(\alpha_j - \alpha_{j+1}). 
\]
\end{lemma}
\begin{proof}
Since $(\alpha_\ell)_{\ell \ge 1}$ is strictly decreasing, $P_\ell>0$ for every $\ell \ge 1$. By expanding the sum, we have $\sum_{\ell \ge i}P_{\ell} =\alpha_i/\alpha_1$, and in particular
$\sum_{\ell \ge 1} P_\ell= 1$. Since $\sum_{\ell \ge 1}\ell P_\ell=\sum_{i \ge 1}\sum_{\ell \ge i}P_{\ell}$ and $\sum_i\alpha_i=1$, we obtain 
\[
p = \frac1{\alpha_1}\sum_{i\ge1}\alpha_i
=\frac1{\alpha_1}<\infty .
\]
Then, from Lemma \ref{RNW}, it follows that $\nu_{\alpha}$ is $\sigma-$invariant, ergodic, and $\nu_{\alpha}(A_i)=\sum_{\ell \ge i}P_{\ell}/p=\alpha_i$ for every $i \ge 1$. Moreover,
\[
\begin{split}
h_{\nu_{\alpha}}(\sigma) = - \frac{1}{p}\sum_{\ell \ge 1}P_\ell \log P_\ell &= - \alpha_{1}\sum_{\ell \ge 1}\frac{\alpha_{\ell} - \alpha_{\ell+1}}{\alpha_1}\log \bigg(\frac{\alpha_{\ell} - \alpha_{\ell+1}}{\alpha_1}\bigg)\\
&=\alpha_1\log \alpha_1 - \sum_{\ell \ge 1}(\alpha_\ell - \alpha_{\ell+1})\log(\alpha_\ell - \alpha_{\ell+1}),
\end{split}
\]
as desired. Finally, ergodicity together with Birkhoff's theorem applied to each $\mathbf 1_{A_j}$ immediately gives $\nu_{\alpha}(C_\alpha)=1$.
\end{proof}

For each $\ell \ge 1$, define the periodic point $z_{\ell}$ with the repeating block $1 0^{\ell-1}$ (the symbol $1$ followed by $\ell -1$ symbols $0$), and consider the probability measure 
\[
\mu_{\ell}:= \frac{1}{\ell}\sum_{j \le \ell}\delta_{\sigma^{j}{(z_{\ell})}}.
\]
The point $z_\ell$ is periodic of least period exactly $\ell$, its orbit consists of $\ell$ distinct points and, since $\mu_{\ell}$ is supported on a finite set, one can easily check that $\mu_\ell$ is a $\sigma$-invariant ergodic probability measure with $h_{\mu_\ell}(\sigma)=0$. Moreover, we also have
\begin{equation}\label{MAJ}
\mu_\ell(A_j)=
\begin{cases}
1/\ell,& 1\le j\le\ell,\\
0,& j>\ell.
\end{cases}
\end{equation}

Recall that $c_0$ is the space of all sequences that converge to zero in the supremum norm. It is well known that the dual space of $c_0$ is $\ell^{1}$.  
\begin{lemma}\label{C0F}
Let $q \in c_0$, consider $f:= \sum_{j \ge 1}q_j\mathbf{1}_{A_j}$, and denote $Q_{\ell}:= \sum_{j \le \ell}q_j$. Then
\begin{itemize}
\item $\displaystyle\int_{\Sigma^{\N}}f\, d\mu_{\ell} = \dfrac{Q_{\ell}}{\ell}$ for every $\ell \ge 1$. 

\item For each measure $\nu$ produced by Lemma \ref{RNW}, we have $\displaystyle\int_{\Sigma^{\N}}f \,d\nu = \frac{1}{p}\sum_{\ell \ge 1}P_{\ell}Q_{\ell}$. 
\end{itemize}
\end{lemma}
\begin{proof}
The first item follows directly by the definition of $f$ together with the definition of $\mu_{\ell}$ and using \eqref{MAJ}, that is, 
\[
\int_{\Sigma^{\N}} f\,d\mu_\ell = \int_{\Sigma^{\N}}\sum_{j \ge 1}q_j\mathbf{1}_{A_j}\, d\mu_{\ell} = \sum_{j \ge 1}q_{j}\mu_{\ell}(A_{j}) 
=\sum_{j\le\ell}\frac{q_j}{\ell}=\frac{Q_\ell}{\ell}. 
\]

For the second item, recalling that $\nu(A_j) = \sum_{\ell \ge j}P_{\ell}/p$ for every $j \ge 1$ (Lemma~\ref{RNW}) and proceeding as before, we can calculate
\[
\begin{split}
\int_{\Sigma^{\N}} f\,d\nu = \int_{\Sigma^{\N}}\sum_{j \ge 1}q_j\mathbf{1}_{A_j}\, d\nu &= \sum_{j \ge 1}q_{j}\nu(A_{j}) = \frac{1}{p}\sum_{j\ge1}q_j\sum_{\ell\ge j}P_{\ell} \\
& = \frac{1}{p}\sum_{\ell\ge1}P_{\ell}\sum_{j\le\ell}q_j = \frac{1}{p}\sum_{\ell\ge1}P_{\ell}Q_{\ell},
\end{split}
\]
and the lemma is proved. 
\end{proof}

The following is the well known \emph{Kac lemma} (see for example Chapter 1 in \cite{VO16}), which we state here for convenience. Given any measure $\eta$, and any subset $S \subset \Sigma^{\N}$ with $\eta(S) > 0$, we denote by $\eta_{S}$ the restriction of $\eta$ to the measurable subsets of $S$. 

\begin{lemma}[Kac formula]
Consider $B:= A_1$ as before, and let $\eta$ be an ergodic $\sigma$-invariant probability measure on
$\Sigma^{\mathbb{N}}$ with $\eta(B)>0$. Then, $\eta_{B}$ is a $\sigma_B-$invariant ergodic
probability measure on $B$ satisfying
\[
\int_{B}r\,d\eta_{B}=\frac{1}{\eta(B)}<\infty.
\]
Moreover, for every bounded measurable function $\psi\colon\Sigma^{\mathbb{N}}\to\mathbb{R}$,
\[
\int_{\Sigma^{\mathbb{N}}}\psi\,d\eta
=\eta(B)\int_{B}\bar\psi\,d\eta_{B},
\qquad\text{where}\qquad
\bar\psi(x):=\sum_{i=0}^{r(x)-1}\psi\bigl(\sigma^{i}(x)\bigr).
\]
\end{lemma}

For a fixed $N \ge 1$ and $q \in \R^{N}$, we denote $\widetilde{q}_j = q_j$ for $j \le N$ and $\widetilde{q}_j = 0$ for $j > N$. Furthermore, we consider 
\[
\xi_{q} := \sum_{j \le N}q_{j}\textbf{1}_{A_j} \quad \textrm{and} \quad Q_k = \sum_{j \le \min\{k,N\}}q_j.   
\]

\begin{lemma}\label{LBD}
There exists a unique $\lambda = \lambda(q) \in (1,\infty)$ such that 
\[
\sum_{k \ge 1}e^{Q_k}\lambda^{-k} = 1 \quad \textrm{and} \quad P_{\rm{top}}(\xi_{q}) = \log \lambda,
\]
and the unique equilibrium state $\nu_q$ of $\xi_{q}$ is a measure produced by Lemma \ref{RNW} with $P_{k} = e^{Q_k}\lambda^{-k}$.  
\end{lemma}
\begin{proof}
For $\lambda>1$, we have that $g(\lambda):= \sum_{k \ge 1} e^{Q_k}\lambda^{-k}$ converges because its tail is the sum $e^{Q_N}\sum_{k\ge N}\lambda^{-k}$. The map $\lambda \mapsto g(\lambda)$ is continuous, strictly decreasing, and it goes to $\infty$ as $\lambda \to 1$. In addition, $g(\lambda)$ tends to $0$ as $\lambda\to\infty$, which readily guarantees the existence of a unique $\lambda$ such that $g(\lambda) = 1$. In this case, $P_k:=e^{Q_k}\lambda^{-k}$ is a probability vector, and we get 
\[
\begin{split}
- \sum_{k \ge 1}P_k\log P_k &= \sum_kP_k(k\log\lambda-Q_k)\le p\log\lambda+ \sup_{k \ge 1}|Q_k| \\
& \le p \log \lambda + (|q_1|+ \cdots + |q_N|) <\infty.
\end{split}
\]

Observe that for $x\in B$ with $r(x)=k$, we have $\sigma^{0}(x) = x \in A_1$ and $\sigma^{i}(x) \in A_{k+1-i}$ for $1\le i\le k-1$. Hence, the orbit segment $x,\sigma(x),\dots,\sigma^{k-1}(x)$ meets each of $A_1,\dots,A_k$ exactly once, and does not meet the other $A_m$. Based on this, one can see that
\[
\bar{\xi}_q(x) = \sum_{i=0}^{r(x)-1}\xi_q\bigl(\sigma^{i}(x)\bigr) = \sum_{i=0}^{r(x)-1}\sum_{j \le N}q_{j}\textbf{1}_{A_j}(\sigma^{i}(x)) = \sum_{j \le k}\widetilde{q}_{j} = Q_{k}. 
\]
Let $\nu$ be the measure produced  via the probability vector $P_k$ above. Now, for the particular measure $\nu_{B} = \bar{\nu}$, we can apply Lemma \ref{RNW} and Kac's formula to obtain 
\begin{equation}\label{EQN1}
\int_{\Sigma^{\N}}\xi_q\,d\nu =\nu(B)\int_{B}Q_{r}\,d\bar{\nu} = \frac1p\sum_{k \ge 1}P_k Q_k.
\end{equation}
Furthermore, we also get 
\begin{equation}\label{EQN2}
h_\nu(\sigma) = \frac{1}{p}\sum_{k \ge 1}P_k(k\log\lambda-Q_k).
\end{equation}
From \eqref{EQN1} and \eqref{EQN2}, we obtain
\begin{equation}\label{EQL} h_\nu(\sigma)+\int_{\Sigma^{\N}}\xi_q\,d\nu=\frac{\log\lambda}{p}\sum_{k \ge 1}k P_k=\log\lambda. 
\end{equation}
In particular, by the variational principle, we have $P_{\rm{top}}(\xi_q) \ge \log \lambda$. 

Let us now show that $\log \lambda$ is an upper bound for the topological pressure of $\xi_q$ and, consequently, that $\nu$ is the equilibrium state of $\xi_q$. In fact, if $\mu \in \cM(\sigma)$ satisfies $\mu(B)=0$, then $\mu(\sigma^{-k}[1])=0$ for all $k$, which forces $\mu=\delta_{0^{\infty}}$, where $0^{\infty}$ is the sequence with all entries equal to 0. In this case, we clearly have 
\[
h_\mu(\sigma) +\int_{\Sigma^{\N}}\xi_q\,d\mu = 0 <\log\lambda.
\]
On the other hand, assume $\mu(B)>0$, and consider $\bar\mu:= \mu_{B}$. Invoking Kac's and Abramov's formulae again, we can see that
\[
h_\mu(\sigma)+\int_{\Sigma^{\N}}\xi_q\,d\mu
=\frac{h_{\bar\mu}(\sigma_B)+\int_B\bar{\xi}_q\,d\bar\mu}{\int_Br\,d\bar\mu} \le \log \lambda 
\]
is equivalent to $h_{\bar\mu}(\sigma_B)+\int_B(\bar{\xi}_q-r\log\lambda)\,d\bar\mu \le 0$. Moreover, for each $x$ with $r(x) = k$, we have $\bar{\xi}_q(x)-r(x)\log\lambda = Q_k-k\log\lambda$. Writing $p_k:=\bar\mu\{r=k\}$, we have $\sum_{k \ge 1} k p_k=\int_{B}r\,d\bar\mu<\infty$. Considering the generating partition $\mathcal{P}_{r}$ in \eqref{PPR}, we have $H_{\bar\mu}(\mathcal{P}_r) = -\sum_{k \ge 1} p_k\log p_k <\infty$ and, by Sinai's Generator Theorem again, $h_{\bar\mu}(\sigma_B) = h_{\bar{\mu}}(\sigma_B,\mathcal{P}_r) \le H_{\bar\mu}(\mathcal{P}_r)$. Finally, taking $c_k:= Q_k-k\log\lambda$, we apply Jensen's inequality to obtain
\[
H_{\bar\mu}(\mathcal{P}_r)+\sum_{k \ge 1}p_kc_k=\sum_{k \ge 1}p_k\log\frac{e^{c_k}}{p_k}
\le\log\sum_{k \ge 1}e^{c_k}=\log\sum_{k \ge 1}e^{Q_k}\lambda^{-k}=0.
\]
That is, 
\[
\begin{split}
h_{\bar\mu}(\sigma_B)\le H_{\bar\mu}(\mathcal{P}_r) &\le -\sum_{k \ge 1}p_kc_k = - \sum_{k \ge 1}p_kQ_k + \log \lambda\sum_{k \ge 1}kp_k \\
& = -\int_{B}\bar{\xi}_q\,d\bar{\mu} + \log \lambda\int_{B}r \,d\bar{\mu}.
\end{split}
\]
Since $\mu$ is arbitrary, we conclude that $P_{\rm{top}}(\xi_q) \le \log \lambda$, as desired. 

Therefore, $P_{\rm{top}}(\xi_q)=\log\lambda$, and \eqref{EQL} guarantees that $\nu$ is an equilibrium state of $\xi_q$. Since $\xi_q$ is locally constant, it has a unique equilibrium state (Theorem 1.22 in \cite{Bow75}), so $\nu = \nu_q$ and the lemma is proved. 
\end{proof}

We are now in a position to finally prove the main result of the section. \\

\emph{Proof of Theorem \ref{NTW}}. We start showing that $\alpha$ belongs to the pseudo interior of $\rm{Rot}(\Xi)$. Fix $N \ge 1$, $w = (w_1,\dots,w_N) \in \R^{N}\setminus\{0\}$, and consider the vector $q:= (w_1,\dots,w_N,0,0,\dots)$. It follows from Lemma \ref{C0F} that
\[
\langle w, \alpha^{(N)}\rangle = \sum_{j \ge 1}q_j\alpha_j =  \frac{1}{p}\sum_{j\ge 1}q_j\sum_{\ell \ge j}P_{\ell} = \frac{1}{p}\sum_{\ell\ge 1}P_{\ell}Q_{\ell} = \sum_{\ell \ge 1}\frac{\ell P_{\ell}}{p}\frac{Q_{\ell}}{\ell}.
\]
Here $P_{\ell}$ and $p$ are as in Lemma \ref{ENT}. Since $Q_{\ell}/\ell$ is bounded and not a constant (independent of $\ell$), we get 
\begin{equation}\label{WAN}
\inf_{\ell \ge 1}\frac{Q_{\ell}}{\ell} = \sum_{\ell \ge 1}\frac{\ell P_{\ell}}{p}\inf_{\ell \ge 1}\frac{Q_{\ell}}{\ell}  < \langle w, \alpha^{(N)}\rangle < \sum_{\ell \ge 1}\frac{\ell P_{\ell}}{p}\sup_{\ell \ge 1}\frac{Q_{\ell}}{\ell} = \sup_{\ell \ge 1}\frac{Q_{\ell}}{\ell}.
\end{equation}
Applying Lemma \ref{C0F} once more, we have
\[
\frac{Q_{\ell}}{\ell} = \int_{\Sigma^{\N}}f\, d\mu_{\ell} = \langle w, (\mu_{\ell}(A_1),\dots, \mu_{\ell}(A_N))\rangle.
\]
As the vector $(\mu_{\ell}(A_1),\dots,\mu_{\ell}(A_N)) \in \rm{Rot}$$(\xi^{(N)})$, it follows from \eqref{WAN} that
\begin{equation}\label{ISQ}
\inf_{y \in \textrm{Rot}(\xi^{(N)})}\langle w, y \rangle \le \inf_{\ell \ge 1}\frac{Q_{\ell}}{\ell} < \langle w, \alpha^{(N)}\rangle < \sup_{\ell \ge 1}\frac{Q_{\ell}}{\ell} \le \sup_{y \in \textrm{Rot}(\xi^{(N)})}\langle w, y \rangle.
\end{equation}
Suppose $\alpha^{(N)} \in \partial \rm{Rot}$$(\xi^{(N)})$. The convexity of $\rm{Rot}$$(\xi^{(N)}) \subset \R^{N}$ guarantees the existence of $w' \neq 0$ such that $\langle w', y \rangle \le \langle w', \alpha^{(N)} \rangle$ for all $y \in \rm{Rot}$$(\xi^{(N)})$, which contradicts \eqref{ISQ}. The same argument also shows that the minimum cannot be attained by $\alpha^{(N)}$. By the arbitrariness of $N$, we conclude that $\alpha^{(N)}$ belongs to $\interior \rm{Rot}$$(\xi^{(N)})$ for every $N \ge 1$, as desired.

We now show that $\alpha \notin \cW$. Fix $N \ge 4$, let $\lambda$ and $P_k=e^{Q_k}\lambda^{-k}$ be as in Lemma~\ref{LBD}. Let $\nu_{q}$ be the unique equilibrium state of $\xi^{(N)}_{q}$, which here is also the measure produced by Lemma \ref{RNW} with $P_k$ as above. For $1\le i\le N-1$, we have
\[
\alpha_i-\alpha_{i+1}=\nu_q(A_i)-\nu_q(A_{i+1})=\frac{P_i}{p}
=\frac{e^{Q_i}\lambda^{-i}}{p}.
\]
Since $q_i= Q_i-Q_{i-1}$, we get
\[
q_i=\log\bigg(\frac{\alpha_i-\alpha_{i+1}}{\alpha_{i-1}-\alpha_i}\bigg)+\log\lambda.
\]
Write $\beta_i:=4^{-2^{i}}$, so $\alpha_i=c\beta_i$ and
\[
\frac{\alpha_i-\alpha_{i+1}}{\alpha_{i-1}-\alpha_i}
=\frac{\beta_i}{\beta_{i-1}}\cdot
\frac{1-\beta_{i+1}/\beta_i}{1-\beta_i/\beta_{i-1}} .
\]
Since $\beta_i/\beta_{i-1}=4^{-2^{i-1}}\le\frac1{16}$ for $i\ge2$ and
$\beta_{i+1}/\beta_i=4^{-2^{i}}\le\frac1{16}$ for $i\ge1$, the second factor
lies in $[\frac{15}{16},\frac{16}{15}]$, while
$\log(\beta_i/\beta_{i-1})=-2^{i-1}\log4$. Hence
\begin{equation}\label{eq:qasym}
q_i=-2^{i-1}\log4+\log\lambda+\delta_i,
\qquad|\delta_i|\le\log\dfrac{16}{15} \quad \textrm{for all} \quad 2\le i\le N-1.
\end{equation}
Moreover $\log\lambda= P_{\rm{top}}(\xi_q)\le h_{\rm{top}}(\sigma)+\|\xi_q\|_\infty
=\log2+\|q\|_\infty$. Taking $i=N-1$ in \eqref{eq:qasym} and using $q_{N-1}\ge-\|q\|_\infty$, we get
\[
-\|q\|_\infty\ \le\ q_{N-1}\ \le\ -2^{N-2}\log4+\log2+\|q\|_\infty + \log\dfrac{16}{15},
\]
which readily implies
\[
\|q\|_\infty\ge2^{N-3}\log4-\dfrac12\bigg(\log2+ \log\dfrac{16}{15}\bigg).
\]
This means that for $N \ge 4$, the vector $\|q^{(N)}\|_{\infty} \ge 2^{N-3} \to \infty$, and $\alpha$ cannot be in $\cW(R)$ for any $R > 0$, that is, $\alpha$ is not uniformly attainable. 

Now let us prove the second item of the theorem. By Lemma \ref{ENT}, we already know that $\nu_{\alpha} \in \cM(\alpha)$, and $\nu_{\alpha}(C_{\alpha}) = 1$. We need to show that $h_{\rm{top}}(C_{\alpha}) = h_{\nu_{\alpha}}(\sigma)$, and that $C_{\alpha}$ is a proper subset of $\Sigma^{\N}$.  

Consider the countable partition
$\Omega:=\{A_1,A_2,\dots\}\cup\{\{0^{\infty}\}\}$ of $\Sigma^{\N}$. If $\mu\in\mathcal M(\alpha)$ then $\sum_{j \ge 1}\mu(A_j)=\sum_{j \ge 1}\alpha_j=1$, which gives $\mu(\{0^{\infty}\})=0$. Moreover, for each $\mu \in \cM(\alpha)$, it is clear that $\Omega$ is a generating partition modulo $\mu$, with entropy 
\[
H_{\mu}(\Omega) = - \sum_{j \ge 1}\alpha_j \log \alpha_j = \log 4 \sum_{j \ge 1}\alpha_{j}2^{j} - \log c \sum_{j \ge 1}\alpha_j < \infty.
\]
Once more we invoke Sinai's Generator Theorem to get $h_{\mu}(\sigma) = h_{\mu}(\sigma,\Omega)$ for each $\mu \in \cM(\alpha)$. Also, observe that, for $\mu \in \cM(\alpha)$, we have 
\[
\mu\bigl(A_j\cap\sigma^{-1}(A_i)\bigr)=
\begin{cases}
\alpha_{i+1},&j=i+1,\\
\alpha_i-\alpha_{i+1},&j=1,\\
0,&\text{otherwise.}
\end{cases}
\]
In this case, considering the partition $\Omega^{2}$ formed by the intersection of elements of $\Omega$ with elements of $\sigma^{-1}\Omega$, we get 
\[
\begin{split}
H_{\mu}(\Omega^{2}) &= -\sum_{i,j \ge 1}\mu(A_j \cap \sigma^{-1}A_i)\log\mu(A_j \cap \sigma^{-1}A_i)\\
&= -\sum_{i\ge 1}\alpha_{i+1}\log \alpha_{i+1} - \sum_{i \ge 1}(\alpha_i - \alpha_{i+1})\log(\alpha_i - \alpha_{i+1})\\
& = H_{\mu}(\Omega) + \alpha_1\log \alpha_1 - \sum_{i \ge 1}(\alpha_i - \alpha_{i+1})\log(\alpha_i - \alpha_{i+1}).
\end{split}
\]
Now we can compute the conditional entropy of $\Omega$ given $\sigma^{-1}\Omega$ 
\[
\begin{split}
H_{\mu}(\Omega | \sigma^{-1}\Omega) &= H_{\mu}(\Omega^{1}) - H_{\mu}(\sigma^{-1}\Omega) = H_{\mu}(\Omega^{1}) - H_{\mu}(\Omega)\\
& = \alpha_1\log \alpha_1 - \sum_{i \ge 1}(\alpha_i - \alpha_{i+1})\log(\alpha_i - \alpha_{i+1}). 
\end{split}
\]
By Lemma \ref{ENT}, this readily implies $H_{\mu}(\Omega | \sigma^{-1}\Omega) = h_{\nu_{\alpha}}(\sigma)$. In particular, for every $\mu\in\mathcal M(\alpha)$, the definition of metric entropy gives
\[
h_\mu(\sigma)\ \le H_{\mu}(\Omega | \sigma^{-1}\Omega)  = h_{\nu_{\alpha}}(\sigma).
\]
This, together with Theorem \ref{ACVP1}, yields 
\[
h_{\rm{top}}(C_{\alpha}) =\sup_{\mu\in\mathcal M(\alpha)}h_\mu(\sigma)
=h_{\nu_{\alpha}}(\sigma),
\]
which proves that $\nu_{\alpha}$ is an ergodic full measure, as desired. Finally, since $\nu_{\alpha}$ is not the Bernoulli measure with two equal weights $1/2$, it cannot be the measure of maximal entropy, that is, $h_{\nu_{\alpha}}(\sigma) < h_{\rm{top}}(\sigma) = \log 2$. Furthermore, by the explicit formula given in Lemma \ref{ENT}, one can check that $h_{\nu_{\alpha}}(\sigma) \approx 0.221$, and the topological entropy of $C_{\alpha}$ is strictly positive. This concludes the proof of the theorem. \qed 

\subsection{Regularity of the entropy spectrum}\label{RES}

Compared to the finite-dimensional setup, the entropy spectrum for the case with countable observables behaves quite differently, as we shall see now. 

Recall that the entropy spectrum is defined as $\cE(\alpha) = h_{\rm{top}}(C_{\alpha})$. Our dynamical system is the shift $\sigma$ on $\{0,1\}^{\N}$. Again, we denote $\Sigma = \{0,1\}$. For a finite word $w$ with length $|w|$, we write $[w]$ for the corresponding cylinder, and let $\mathcal{V}$ be the set of all finite words of length at least one. Since $\mathcal{V}$ is countable, we can fix any enumeration $\mathcal{V} = \{w_1, w_2, \dots\}$. Now we consider the observable $\Phi \colon \Sigma^{\N} \to \R^{\N}$ with coordinates given by 
\[
\phi_{j}:= 8^{-|w_j|}\textbf{1}_{[w_j]}. 
\]

The main result of this subsection is the following:

\begin{theorem} \label{REG}
For the observable $\Phi$, the following assertions hold: 
\begin{enumerate}
\item $\cE$ is affine and upper semicontinuous on $\rm{Rot}(\Phi)$ (wrt the product topology), and $0 \le \cE \le \log 2$, where the maximum and minimum are attained. 

\item $\{\alpha \in \rm{Rot}(\Phi): \cE(\alpha) = 0\}$ is a dense $G_{\delta}$ subset and, $\cE$ is continuous at each of its points. 

\item $\cE$ is discontinuous at every $\alpha$ with $\cE(\alpha) > 0$. 

\item For each $\alpha_0, \alpha_1 \in \rm{Rot}(\Phi)$ the function $t \mapsto \cE((1-t)\alpha_0 + t\alpha_1)$ is continuous on $[0,1]$. 
\end{enumerate}
\end{theorem} 

\begin{proof}
Since there are $2^m$ words of length $m$,
\[
\sum_{j\ge1}8^{-|w_j|}=\sum_{m\ge1}2^m8^{-m}=\frac13 ,
\]
so $\Phi$ takes values in the Banach space $\ell^1$. By dominated convergence, coordinatewise convergence of $\frac{1}{N}S_N\Phi(x)$ is the same as convergence in $\ell^1$. In this case, it does not matter which of the two ways is used to define the level set $C_\alpha$. Now consider the map $\Gamma \colon\mathcal M(\sigma)\to\ell^1$ given by 
\[
\Gamma(\mu):= \Biggl(\int_{\Sigma^{\N}}\phi_1\,d\mu, \int_{\Sigma^{\N}}\phi_2\,d\mu, \dots \Biggr) = \bigg(8^{-|w_1|}\mu([w_1]), 8^{-|w_2|}\mu([w_2]), \dots\bigg). 
\]
Observe that $\mathrm{Rot}(\Phi)=\Gamma(\mathcal M(\sigma))$. The map $\Gamma$ is affine, and it is injective because two invariant measures agreeing on every cylinder are equal. It is also easy to check that $\Gamma$ is continuous from the weak$^*$ topology to the norm of $\ell^1$. Hence, $\Gamma$ is an affine homeomorphism from $\mathcal M(\sigma)$ onto $\mathrm{Rot}(\Phi)$, and $\mathrm{Rot}(\Phi)$ is compact and convex. For $\alpha\in\mathrm{Rot}(\Phi)$ we denote by $\mu_\alpha:=\Gamma^{-1}(\alpha)$ the unique invariant measure with $\int_{\Sigma^{\N}}\Phi\,d\mu_\alpha=\alpha$.

Let us now identify the level sets. Consider the empirical measures of $x \in \Sigma^{\N}$
\[
\nu_N(x):=\frac{1}{N}\sum_{k = 0}^{N-1}\delta_{\sigma^{k}(x)},
\]
and $V(x)$ the nonempty compact set of their limit
points with respect to the weak$^{*}$ topology. Recall that $G_\mu:=\{x \in \Sigma^{\N}:\ \nu_N(x)\to\mu\}$ is the set of generic points of $\mu\in\mathcal{M}(\sigma)$. Since
$\int_{\Sigma^{\N}}\phi_j\,d\nu_N(x)=\frac{1}{N}S_N\phi_j(x)$ and each $\phi_j$ is
continuous, every $\mu\in V(x)$ with $x\in C_\alpha$ satisfies
$\Gamma (\mu)=\alpha$, which implies $\mu=\mu_\alpha$. Hence, $V(x)=\{\mu_\alpha\}$ and $x\in G_{\mu_\alpha}$. Conversely, $\nu_N(x)\to\mu_\alpha$ gives $\frac{1}{N}S_N\phi_j(x)\to\alpha_j$ for every $j \ge 1$. Therefore
$C_\alpha=G_{\mu_\alpha}$ for every $\alpha\in\mathrm{Rot}(\Phi)$. The full shift has the specification property, so it is saturated, that is, $h_{\mathrm{top}}(\sigma|G_\mu)=h_\mu(\sigma)$ for every $\mu\in\mathcal{M}(\sigma)$ (\cite{FLP08,PS07}). We then conclude that
\begin{equation}\label{eq:spectrum}
\mathcal{E}(\alpha) = h_{\rm{top}}(C_{\alpha}) = h_{\mathrm{top}}(\sigma|G_{\mu_\alpha})=h_{\mu_\alpha}(\sigma),
\end{equation}
that is, $\mathcal{E} = h\circ\Gamma^{-1}$, where $h\colon\mathcal M(\sigma)\to[0,\log2]$ is the entropy map $h(\mu) := h_{\mu}(\sigma)$. All four
assertions in the theorem are now properties of $h$ transported by the affine homeomorphism $\Gamma$. 

Indeed, $h$ is upper semicontinuous and affine on $\cM(\sigma)$. Thus, because of $\Gamma$, the spectrum $\mathcal{E}$ is forced to satisfy the same properties. Moreover $h_\mu(\sigma) \le\frac{1}{N}\log2^N=\log2$, with
equality attained at the Bernoulli measure $m_{1/2}$, with the two weights equal to 1/2. Also note that $h(\delta_{0^\infty})=0$.  Hence, it readily follows from \eqref{eq:spectrum} that $\mathcal{E}$ satisfies $0\le\mathcal E\le\log2$, with the maximum and minimum of the range attained. This proves (1).

For the item (2), note that 
\[
\Big\{\alpha \in \rm{Rot}(\Phi):\ \mathcal{E}(\alpha)=0\Big\}=\Gamma(Z)
\]
with $Z:=\{\mu \in \mathcal{M}(\sigma):\ h(\mu) = 0\}$. By upper semicontinuity, each set $\{h<1/k\}$ is open, and so
$Z=\bigcap_{k \ge 1}\{h<1/k\}$ is a $G_\delta$ set. Furthermore, $Z$ is dense because, for the full shift, it is well known that the invariant measures supported on periodic orbits are dense in $\mathcal M(\sigma)$ and have zero entropy. Since both properties are preserved by the homeomorphism $\Gamma$, we conclude that $\Gamma(Z)$ is a dense $G_\delta$ subset of $\mathrm{Rot}(\Phi)$.

For the item (3), let $\mathcal E(\alpha)>0$. By the density obtained in item (2), there are points $\alpha_N\in\mathrm{Rot}(\Phi)$ with $\alpha_N\to\alpha$ and $\mathcal E(\alpha_N)=0$ for every $N \ge 1$. Thus, $\lim_{N \to \infty}\mathcal {E}(\alpha_N) = 0 < \mathcal E(\alpha)$. 

Finally, we observe that item (4) is a consequence of affinity plus convexity. In fact, since $\mathrm{Rot}(\Phi)$ is
convex, for $t\in[0,1]$ and $\alpha_0, \alpha_1 \in \rm{Rot}(\Phi)$, we have
\[
\mathcal {E}\bigl((1-t)\alpha_0+t\alpha_1\bigr)=(1-t)\mathcal{E}(\alpha_0)+t\,\mathcal E(\alpha_1),
\]
which is, of course, continuous with respect to $t$. 
\end{proof}

\begin{remark}
The only property of $\Phi$ strongly used in the proof of Theorem \ref{REG} is that the integration map $\Gamma\colon \cM(\sigma) \to \R^{\N}$ is injective (the summability property $\sum_j\|\phi_j\|_\infty < \infty$ only serves to identify the $\ell^1$ and product topologies on $\rm{Rot}(\Phi)$). For any such \emph{separating} observable, $\Gamma$ is an affine homeomorphism from $\cM(\sigma)$ onto $\rm{Rot}(\Phi)$, $C_\alpha = G_{\mu_\alpha}$ with $\mu_\alpha := \Gamma^{-1}(\alpha)$, and $\cE = h\circ\Gamma^{-1}$, so Theorem \ref{REG} holds verbatim. In particular, it holds for the observable of Theorem \ref{NEFM}, which is also separating by Lemma~\ref{lem:sep} and has nonempty pseudo-interior. Note also that for a separating observable we have $\cM(\alpha) = \{\mu_\alpha\}$, hence $\alpha$ admits a full measure if and only if $\mu_\alpha$ is ergodic, in which case it is unique. Observe that the parameter $\alpha$ of Theorem \ref{NEFM} is an instance of this dichotomy with $\mu_\alpha = \eta$ non-ergodic.
\end{remark}

\begin{remark}
As noted before, the observable of Theorem \ref{NEFM} may be taken of the form $\phi_j = c_j\,\mathbf 1_{[w_j]}$ for a suitable sparse family of words, with any summable positive weights $c_j$. In contrast, the family of \emph{all} cylinders used in Theorem \ref{REG} is too large for Theorem \ref{NEFM}: the identities $\mathbf 1_{[0]} + \mathbf 1_{[1]} \equiv 1$ and $\mathbf 1_{[w0]} + \mathbf 1_{[w1]} = \mathbf 1_{[w]}$ force $\rm{Rot}$$(\phi^{(N)})$ into an affine hyperplane for every large $N$, so that $\interior_{\rm ps}\rm{Rot}(\Phi) = \emptyset$ for that observable. Items (1), (3), (4) and (5) of Theorem \ref{NEFM} still hold for it, since they only use that the integration map separates invariant measures together with $\sum_j\|\phi_j\|_\infty + \sum_j\rm{Lip}(\phi_j) < \infty$. Item (2) needs the existence of pseudo-interior points and this is exactly the contribution of Lemmas \ref{lem:algebra} and \ref{lem:complement}. 
\end{remark}
	
\section{Some applications of Theorems \ref{ACVP1} and \ref{ACVP2}}\label{ACV}
	
Assuming that $X$ is a metric space with distance $d$, we recall that a function $g\colon X \to \R$ is \emph{Lipschitz continuous} with \emph{Lipschitz constant} $\textrm{Lip}(g) > 0$ if 
\[
\big|g(x) - g(y)\big| \le \textrm{Lip}(g)\,d(x,y) \quad \textrm{for every $x,y \in X$.}
\]
We denote $\textrm{Lip}(X,\R)$ as the space of real-valued Lipschitz continuous functions on $X$. The following simple result is useful for many grounded applications of Theorem~\ref{ACVP2}.
	
\begin{proposition}\label{prop:ex1}
Let $(X,T)$ be a permissible system with $\cU=\mathrm{Lip}(X,\R)$, let $R>0$ and
suppose there are parameters $\lambda_{j}\ge R$ with
\[
\sum_{j\ge1}\lambda_{j}\|\phi_{j}\|_{\infty}<\infty
\qquad\text{and}\qquad
\sum_{j\ge1}\lambda_{j}\,\mathrm{Lip}(\phi_{j})<\infty .
\]
Then every $\alpha\in\mathcal \cW(R)$ admits a unique ergodic full measure.
\end{proposition}
	
\begin{proof}
Since $\alpha_{j}=\int\phi_{j}\,d\mu$ for some $\mu\in \cM(T)$, one has $|\alpha_{j}|\le\|\phi_{j}\|_{\infty}$ and hence $\|\phi_{j}-\alpha_{j}\|_{\infty}\le2\|\phi_{j}\|_{\infty}$. Therefore, for $\gamma\in Q_{R}$, we get
\[
\Bigl\|\sum_{j>N}\gamma_{j}(\phi_{j}-\alpha_{j})\Bigr\|_{\infty}
\le R\sum_{j>N}\|\phi_{j}-\alpha_{j}\|_{\infty}
\le2\sum_{j>N}\lambda_{j}\|\phi_{j}\|_{\infty}\longrightarrow0,
\]
where we recall $Q_R = \{\gamma \in \R^{\N}: \sup_{j \ge 1}|\gamma_j|\le R\}$. The arbitrariness of $\gamma$ yields (A2). In addition, setting $\psi_{\gamma} := \sum_{j \ge 1}\gamma_j(\phi_j - \alpha_j)$, it is clear that $\mathrm{Lip}(\psi_{\gamma})\le\sum_{j}\lambda_{j}\mathrm{Lip}(\phi_{j})<\infty$, and then every $\psi_{\gamma}$ with $\gamma\in Q_{R}$ is Lipschitz and has a unique equilibrium state, which is (A3). The existence of a unique ergodic full measure for $\alpha$ is now guaranteed by Theorem~\ref{ACVP2}. 
\end{proof}
	
\begin{remark}
The class of observables in Proposition \ref{prop:ex1} takes values in $\ell^{1}$. It is then natural to wonder if, in the case of permissible systems, we can always guarantee the existence of full measures for observables taking values in other $\ell^{p}$ spaces, or in other Banach spaces in general. 
\end{remark}
	
To better grasp the quantitative role of uniformly attainable vectors, we exhibit a concrete illustration on the full shift.
	
\begin{example}\label{CPTH}
Let $\Sigma=\{0,1\}$, and consider the shift map $\sigma\colon \Sigma^{\N} \to \Sigma^{\N}$. We recall that for each finite word $\omega = (\omega_1,\dots,\omega_{\ell})$, we denote by $[\omega]$ the cylinder generated by $\omega$. Now let $\Phi\colon \Sigma^{\N} \to \R^{\N}$ be the observable composed of the coordinate functions 
\[
\phi_{j}:=4^{-j}\,\mathbf 1_{[0^{j-1}1]} \quad \textrm{for every $j\ge1$}.
\]
In this case, we immediately get $\|\phi_{j}\|_{\infty}=4^{-j}$ and $\mathrm{Lip}(\phi_{j})\le4^{-j}2^{j}
=2^{-j}$. Thus, Proposition~\ref{prop:ex1} holds with each $\lambda_{j}\equiv R$ for any $R>0$. This happens because $\sum_{j\ge 1}R\,4^{-j}=R/3$ and $\sum_{j\ge 1}R\,2^{-j}=R$.
		
Now take any $t\in[-1,1]$ and let $\nu^{*}(t)$ be the unique equilibrium state of the function $t\phi_{1}=(t/4)\mathbf 1_{[1]}$, that is, the Bernoulli measure with weights
\[
p_{1}=\frac{e^{t/4}}{1+e^{t/4}} \quad \textrm{and} \quad p_{0}=\frac{1}{1+e^{t/4}}.
\]
Take the vector $\alpha \in \R^{\N}$ with coordinates $\alpha_{j} := \alpha_j(t) = \int_{\Sigma^{\N}}\phi_{j}\,d\nu^{*}(t)=4^{-j}p_{0}^{\,j-1}p_{1}$. We claim that $\alpha\in\mathcal \cW(1)$. In fact, for each $N\ge1$ the choice $q^{(N)}=(t,0,\dots,0)$ gives $\sum_{j\le N}q^{(N)}_{j}\phi_{j}=t\phi_{1}$, whose equilibrium state is $\nu^{*}(t)$, and its first $N$ barycentres are exactly $(\alpha_{1},\dots,\alpha_{N})$. Note that by Proposition \ref{OPT} we also have $\alpha \in \interior_{\rm{ps}}\rm{Rot}(\Phi)$. Now it follows from Proposition~\ref{prop:ex1} that $\alpha$ admits a unique ergodic full measure, which is the equilibrium state of $\sum_{j\ge 1}q^{(N)}_{j}\phi_j = t\phi_1$. Thus, in this particular case, the ergodic full measure is exactly $\nu^{*}(t)$, and the entropy spectrum at $\alpha$ is 
\[
h_{\rm{top}}\bigl(C_{\alpha}\bigr)=h_{\nu^{*}}(\sigma)
=-p_{0}\log p_{0}-p_{1}\log p_{1}.
\]
For instance, taking $t=1$ gives $p_{1} \approx 0.562$, $p_{0}\approx 0.438$, and the vector $\alpha \approx (0.141,0.015,\dots) \in \interior_{\rm{ps}}\rm{Rot}(\Phi)$ with entropy $h_{\rm{top}}(C_{\alpha}) \approx 0.685 <\log2 = h_{\rm{top}}(\sigma)$, which means the level set is non-degenerate and proper. 

Observe that, by the construction above, when $t$ varies in $[-1,1]$ it creates a connected path $\nu^{*}(t)$ in the set of shift invariant measures, which creates another connected path $\alpha(t):= (\alpha_1(t),\alpha_2(t),\dots)$ in $\cW(1)$. Since each $\nu^*(t)$ is Bernoulli, the spectrum can be computed explicitly on the curve $\alpha(t)$: 
\[
\cE(t):= h_{\rm{top}}(C_{\alpha(t)}) = \log\bigl(1+e^{t/4}\bigr)-\frac{te^{t/4}}{4(1+e^{t/4})}.
\]
After some sort of ``renormalization'', it behaves like a one-dimensional spectrum for a typical locally constant observable.  
\begin{figure}[ht!]
\centering
\begin{tikzpicture}
\begin{axis}[
width=7.7cm, height=5.04cm,
xlabel={$t$},
ylabel={$\cE(t)$},
ylabel style={rotate=-90},
xmin=-1.05, xmax=1.05,
ymin=0.66, ymax=0.6960,
xtick={-1,-0.75,-0.5,-0.25,0,0.25,0.5,0.75,1},
ytick={0.66,0.67,0.68,0.693147},
yticklabels={$0.66$,$0.67$,$0.68$,$\log 2$},
scaled y ticks=false,
yticklabel style={/pgf/number format/fixed, /pgf/number format/precision=3}, grid=major, grid style={black!20, very thin}, axis line style={black}, tick style={black}, label style={black}, tick label style={black, font=\scriptsize}, label style={font=\small}, x tick label style={/pgf/number format/fixed},]

\addplot[
color=black,
line width=1.6pt, dashed, dash pattern=on 5pt off 3pt, domain=-1:1, samples=2,] {0.693147};

\addplot[
color=black,
line width=2pt, smooth,] coordinates 
{
(-1.0000,0.685395) (-0.9750,0.685775) (-0.9500,0.686146) (-0.9250,0.686507) (-0.9000,0.686859) (-0.8750,0.687201)
(-0.8500,0.687534) (-0.8250,0.687858) (-0.8000,0.688172) (-0.7750,0.688477) (-0.7500,0.688772) (-0.7250,0.689058)
(-0.7000,0.689334) (-0.6750,0.689600) (-0.6500,0.689857) (-0.6250,0.690105) (-0.6000,0.690343) (-0.5750,0.690571)
(-0.5500,0.690789) (-0.5250,0.690998) (-0.5000,0.691198) (-0.4750,0.691388) (-0.4500,0.691568) (-0.4250,0.691738)
(-0.4000,0.691899) (-0.3750,0.692050) (-0.3500,0.692191) (-0.3250,0.692323) (-0.3000,0.692445) (-0.2750,0.692557)
(-0.2500,0.692659) (-0.2250,0.692752) (-0.2000,0.692835) (-0.1750,0.692908) (-0.1500,0.692971) (-0.1250,0.693025)
(-0.1000,0.693069) (-0.0750,0.693103) (-0.0500,0.693128) (-0.0250,0.693142) (0.0000,0.693147) (0.0250,0.693142)
(0.0500,0.693128) (0.0750,0.693103) (0.1000,0.693069) (0.1250,0.693025) (0.1500,0.692971) (0.1750,0.692908)
(0.2000,0.692835) (0.2250,0.692752) (0.2500,0.692659) (0.2750,0.692557) (0.3000,0.692445) (0.3250,0.692323)
(0.3500,0.692191) (0.3750,0.692050) (0.4000,0.691899) (0.4250,0.691738) (0.4500,0.691568) (0.4750,0.691388)
(0.5000,0.691198) (0.5250,0.690998) (0.5500,0.690789) (0.5750,0.690571) (0.6000,0.690343) (0.6250,0.690105)
(0.6500,0.689857) (0.6750,0.689600) (0.7000,0.689334) (0.7250,0.689058) (0.7500,0.688772) (0.7750,0.688477)
(0.8000,0.688172) (0.8250,0.687858) (0.8500,0.687534) (0.8750,0.687201) (0.9000,0.686859) (0.9250,0.686507)
(0.9500,0.686146) (0.9750,0.685775) (1.0000,0.685395)
};

\addplot[only marks, mark=*, mark size=2.2pt, color=black]
coordinates {(0,0.693147)};

\end{axis}
\end{tikzpicture}
\caption{The entropy spectrum is analytic on $\alpha(t)$, with maximum $\cE(0)=\log2$ attained at the measure of maximal entropy $\nu^{*}(0)$, and endpoints $\cE(\pm1)\approx 0.686$.}\label{LSPEC}
\end{figure}
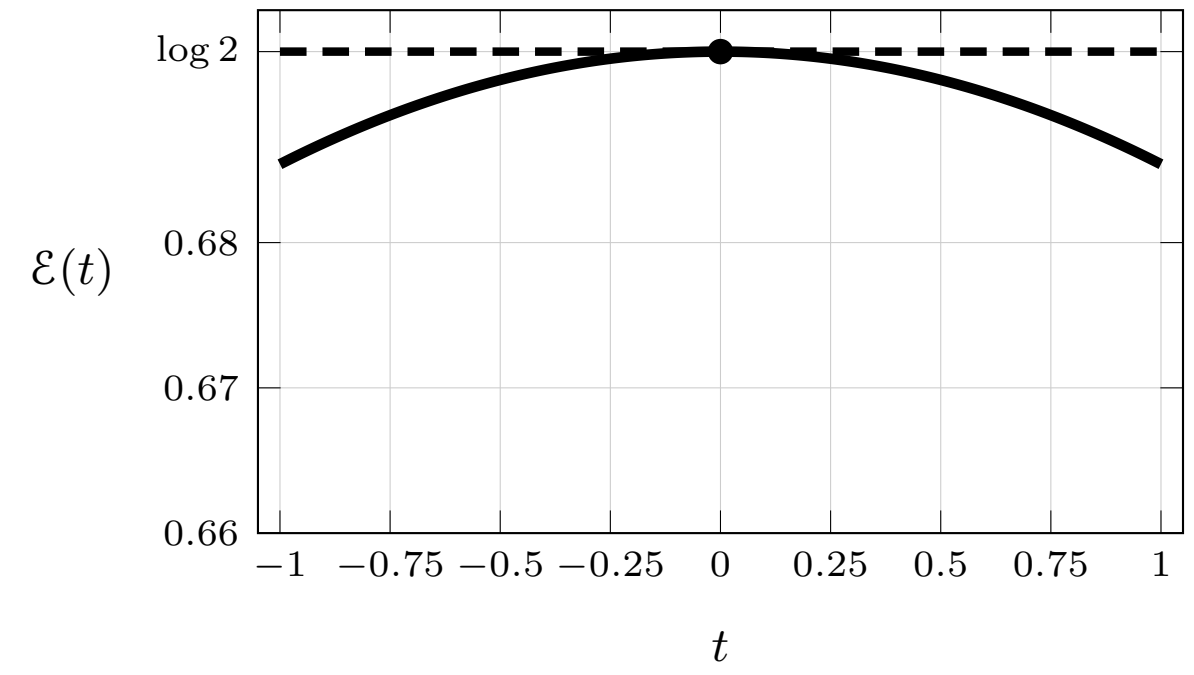   
\end{example}

\subsection{Revisiting recurrence in cylinders: looking for ergodic full measures}
	
This section is inspired by ``recurrence in infinite cylinders'' described by Fan, Liao and Peyri\`ere in \cite{FLP08}. We will show how Theorems \ref{ACVP1} and \ref{ACVP2} can be used to improve information on measures satisfying the conditional variational principle. As discussed earlier, we emphasyze here that the thermodynamic method brings new precise information about the existence of ergodic measures supported on the level sets.  
	
As before, let $\Sigma$ be a finite set and $\sigma: \Sigma^{\N} \to \Sigma^{\N}$ the usual shift map. We recall that the cylinder generated by a word $\omega = (\omega_1,...,\omega_m)$ of length $|\omega| = m$ is the set defined by
\[
C(\omega) := [\omega] = \{\omega' = (\omega'_{1}, \omega'_{2},...) \in \Sigma^{\N}: \omega'_{1} = \omega_1,..., \omega'_{m} = \omega_{m}\}. 
\]
For a given finite word $\omega$, we naturally define the \emph{frequency} of a point $x \in \Sigma^{\N}$ in the cylinder $C(\omega)$ by the limit
\[
\cF(x,C(\omega)):= \lim_{N \to \infty}\frac{1}{N}\sum_{k=0}^{N-1}\textbf{1}_{C(\omega)}(\sigma^{k}(x)) = \lim_{N \to \infty}\frac{1}{N}S_N\textbf{1}_{C(\omega)}(x), 
\]
when it exists.
Let $\cC = \{C(\omega^{j})\}_{j \in \N}$ be a countable collection of cylinders generated by different words $\omega^{j}$ of any length. For a given sequence $\alpha = (\alpha_1,\alpha_2,...)$ with $\alpha_j \in [0,1]$ for each $j \ge 1$, we consider the level set
\[
\cF(\alpha,\cC) := \{x \in \Sigma^{\N}: \cF(x,C(\omega^{j})) = \alpha_j \textrm{ for all $j \ge 1$}\}. 
\]
This set puts together all the points with certain prescribed frequencies. 
	
\subsubsection{General recurrence and approximations by ergodic measures}
	
Fixing a collection of cylinders $\cC$, we consider the function $\textbf{1}_{\cC}: \Sigma^{\N} \to \ell^{\infty}$ given by 
\[
\textbf{1}_{\cC}(x):= (\textbf{1}_{C^1}(x), \textbf{1}_{C^2}(x),\dots,\textbf{1}_{C^n}(x),\dots),
\] 
where, for convenience, we denote by $C^{j}$ the distinct cylinders generated by distinct words $\omega^{j}$. In this case, the rotation set is  
\[
\textrm{Rot}(\textbf{1}_{\cC}) = \bigg\{\bigg(\int_{\Sigma^{\N}}\textbf{1}_{C^{1}}\,d\mu,\int_{\Sigma^{\N}}\textbf{1}_{C^{2}}\,d\mu,\int_{\Sigma^{\N}}\textbf{1}_{C^{3}}\,d\mu, \dots \bigg): \mu \in \cM(\sigma)\bigg\} \subset [0,1]^\N.
\]
The system $(\Sigma^{\N},\sigma)$ is permissible, and it is not hard to see that Theorem \ref{ACVP1}, together with its proof, yields the following result.
	
\begin{corollary}\label{BEF}
Let $\alpha \in \interior_{\rm{ps}}\rm{Rot}(\mathbf{1}_{\cC})$. Then, for each $\delta >0$ there exists an ergodic measure $\eta_{\alpha} \in \cM(\sigma)$ such that $|h_{\rm{top}}(\cF(\alpha, \cC)) - h_{\eta_{\alpha}}(\sigma)| < \delta$. Moreover, there exist $N$ and a vector $q = (q_1,\dots,q_N)$ such that $\eta_{\alpha}$ is the unique equilibrium state for the observable $\sum_{j \le N}q_j(\mathbf{1}_{C^j}-\alpha_j)$. 
\end{corollary}
	
\begin{remark}
This corollary complements the results of Section 5 in \cite{FLP08}. In that paper, the authors noted that the measure of maximal entropy for the entropy spectrum is not ergodic in general. On the other hand, Corollary~\ref{BEF} provides more information and guarantees that for all parameters $\alpha$ in the pseudo-interior of the rotation set, the measure of maximal entropy is nearly supported on the level set of frequencies and can be well approximated by ergodic measures, which are actually unique equilibrium states for some appropriate real-valued observables. 
\end{remark}
	
The natural follow-up question is: can we find an appropriate type of recurrence in infinite cylinders where one can guarantee the existence of an ergodic full measure for some parameter $\alpha$ in the rotation set? The answer is yes, and it happens for certain collections of cylinders generated by words with some specific prescribed lengths, as demonstrated in the next section.
	
\subsubsection{Uniformly normalized recurrence and the existence of ergodic full measures}\label{NRC}

Let $f \colon \N \to \N$ be a non-constant function, and consider a countable collection of finite words $\{\omega^1, \omega^2,\omega^3,\dots\}$ with $|\omega^\ell| = f(\ell)$ for each $\ell \ge 1$. Now take the collection of cylinders with respect to $f$ given by $\cC_{f} := (C^{\ell})_{\ell \in \N}$ where each $C^{\ell}$ is generated by the distinct word $\omega^{\ell}$. Fix a number $\tau \ge 2$ and consider the observable $\Theta\colon \Sigma^{\N} \to \ell^{\infty}$ with each coordinate function defined as  
\[
\theta_{\ell}(x):= \frac{\mathbf{1}_{C^{\ell}}(x)}{\tau^{f(\ell) + \ell -1}} \quad \textrm{for all $x \in \Sigma^{\N}$}. 
\]
Compared to the last section, note how $\Theta$ is the observable $\mathbf{1_{\cC}}$ normalized by the lengths of the words. 
    
One can easily see that 
\[
\sum_{\ell = 1}^{\infty}\|\theta_{\ell}\|_{\infty} \le \sum_{\ell = 1}^{\infty}\frac{1}{\tau^{\ell}} \le \sum_{\ell = 1}^{\infty}\frac{1}{2^{\ell}} < \infty. 
\]
Moreover, for each $\ell \ge 1$ we have 
\[
\big|\theta_{\ell}(x) - \theta_{\ell}(y)\big| = \frac{1}{\tau^{f(\ell) + \ell -1}}\big|\mathbf{1}_{C^{\ell}}(x) - \mathbf{1}_{C^{\ell}}(y)\big| \le \frac{2^{f(\ell)-1}}{\tau^{f(\ell) + \ell -1}}d(x,y) \le \frac{1}{2^{\ell}}d(x,y)
\]
for every $x, y \in \Sigma^{\N}$, implying that $\sum_{\ell = 1}^{\infty}\rm{Lip}(\theta_{\ell}) < \infty$. Thus, the observable $\Theta$ satisfies the conditions of Proposition \ref{prop:ex1}. Now it follows from Theorem \ref{ACVP2} that for each $\alpha \in \rm{Rot}(\Theta) \cap \cW(1)$ there exists a unique ergodic measure $\nu_{\alpha} \in \cM(\alpha)$ satisfying $\nu_{\alpha}(\cF(\alpha,\cC_{f})) = 1$ and such that $h_{\rm{top}}(\cF(\alpha,\cC_{f})) = h_{\nu_{\alpha}}(\sigma)$. Moreover, $\nu_{\alpha}$ is the unique equilibrium state for the observable
\[
\sum_{\ell = 1}^{\infty}\frac{q_{\ell}}{\tau^{f(\ell) + \ell -1}}(1_{C^{\ell}}-\alpha_\ell) \quad \textrm{with $q_{\ell} \in [-1,1]$ for every $\ell \ge 1$}.
\]
In other words and in particular, every $\alpha \in \cW(1)$ admits a unique ergodic full measure. Furthermore, observe that $\alpha \in \interior_{\rm{ps}}\rm{Rot}(\Theta)$, that is, 
\[
\alpha_\ell \in \bigg(0, \frac{1}{\tau^{f(\ell) + \ell -1}}\bigg) \quad \textrm{for every $\ell \ge 1$}.
\]
But not every vector in the above set is in $\cW(1)$, as expected. \qed 
\bigskip

In general it is hard to see the set of uniformly attainable vectors. However, for some particular cases of recurrence in cylinders, by using the machinery from Section \ref{OTW}, we can describe $\cW(1)$ explicitly. In fact, in the recurrence model above, let for instance $\Sigma = \{0,1\}$, $f(\ell)=\ell$ and $\omega^{\ell}=0^{\ell-1}1$, so that $C^{\ell}=A_{\ell} = [0^{\ell-1}1]$ and
\[
\theta_\ell=c_\ell\,\mathbf 1_{A_\ell},\qquad c_\ell:=\tau^{-(2\ell-1)}.
\]
Write $\alpha_j=c_ja_j$, so that $a_j=\nu(A_j)$ for any $\nu\in\cM(\alpha)$. Since $A_{j+1}=[0]\cap\sigma^{-1}(A_j)$, for every $\sigma-$invariant measure $\mu$ we get $\mu(A_{j+1})\le\mu(A_j)$, and $\sum_{j \ge 1}\mu(A_j)=1-\mu(\{0^\infty\})\le1$. Now fix $N\ge1$ and $q\in[-1,1]^N$, and set $t_j:=q_jc_j$, so that $\sum_{j\le N}q_j\theta_j=\sum_{j\le N}t_j\mathbf 1_{A_j}=:\xi_t$. By Lemma~\ref{LBD}, the equilibrium state $\nu_t$ of $\xi_t$ is the induced measure of Lemma~\ref{RNW} with $P_k=e^{Q_k}\lambda^{-k}$, where $Q_k=\sum_{j\le\min\{k,N\}}t_j$ and $\lambda>1$ is the root of $\sum_{k\ge1}e^{Q_k}\lambda^{-k}=1$. In particular
\[
a_j:=\nu_t(A_j)=\frac1p\sum_{\ell\ge j}P_\ell,\qquad p=\sum_{k\ge1}kP_k .
\]
Since $Q_k=\sum_{j\le\min\{k,N\}}t_j$ satisfy $Q_k=Q_N$ for every $k\ge N$, we readilly obtain $P_k=e^{Q_k}\lambda^{-k}=e^{Q_N}\lambda^{-k}$, which gives
\[
\frac{P_{k+1}}{P_k}=\frac1\lambda
\quad\text{and}\quad
P_k=P_N\,\lambda^{-(k-N)}\quad\text{for all }k\ge N.
\]
In particular, for $k\ge N$, we get
\[
\sum_{\ell\ge k}P_\ell=P_k\,\frac{\lambda}{\lambda-1},
\]
which gives $a_j=\frac1p\sum_{\ell\ge j}P_\ell=a_N\lambda^{-(j-N)}$ for $j\ge N$. Now using the identity $\sum_{j\ge1}a_j=1$, we can write
\[
\lambda=\lambda_N(a):=1+\frac{a_N}{1-\sum_{j\le N}a_j}.
\]
Moreover $p=1/a_1$, $P_j=(a_j-a_{j+1})/a_1$ for $j<N$, and $P_N=a_N(1-\lambda^{-1})/a_1$. Using $t_j=Q_j-Q_{j-1}=\log P_j-\log P_{j-1}+\log\lambda$ we obtain the vector $(t_1,\dots,t_N)$ with respect to $\alpha^{(N)} = (\alpha_1,\dots,\alpha_N)$ as 
\begin{equation}\label{eq:tofa}
\begin{aligned}
&t_1:=t_1(a^{(N)})=\log\lambda_N(a)+\log\Big(1-\frac{a_2}{a_1}\Big),\\
&t_i:= t_i(a^{(N)})=\log\lambda_N(a)+\log\bigg(\frac{a_i-a_{i+1}}{a_{i-1}-a_i}\bigg)\qquad \textrm{for $2\le i\le N-1$},\\
&t_{N}:= t_N(a^{(N)})=\log\Bigg(\frac{a_N^{2}}{(a_{N-1}-a_N)\bigl(1-\sum_{j\le N}a_j\bigr)}\Bigg) .
\end{aligned}
\end{equation}
Conversely, given $a_1>a_2>\dots>a_N>0$ with $\sum_{j\le N}a_j<1$, one can see that the vector $t^{(N)} \in \R^{N}$ defined by \eqref{eq:tofa} produces, through Lemma~\ref{LBD}, an equilibrium state with $\nu(A_j)=a_j$ for $j\le N$. Therefore
\[
\cW_N(1)=\Big\{(c_1a_1,\dots, c_Na_N):\ a_1>\dots>a_N>0,\ \ \sum_{j\le N}a_j<1,\ \ |t_i|\le c_i \Big\},
\]
and $\cW(1)$ is then formed of vectors $\alpha = (c_1a_1, c_2a_2, \dots)$ such that $(c_1a_1,\dots, c_Na_N) \in \cW_N(1)$ for all $N \ge 1$. Note that, in accordance with Theorem~\ref{ACVP2}, every $\alpha\in\cW(1)$ admits a unique ergodic full measure, which is the equilibrium state of $\sum_{j\ge1}t^*_j\mathbf 1_{A_j}$ for some $t^{*} = (t_1^{*}, t^{*}_2,\dots)$ with $|t^*_j|\le c_j$ for every $j \ge 1$.

We also note that $\cA_{\rm{D}} = \emptyset$. In fact, $\|\theta_\ell\|_\infty=c_\ell=\tau^{-(2\ell-1)}\to0$, which implies $\sum_\ell\|\theta_\ell\|_\infty<\infty$ and Proposition \ref{OSC} gives $\delta(\alpha)=0$ for every $\alpha$.

\subsection{Realization of vectors as Gibbs states}
	
Based on the formalism explored in the proof of Theorem \ref{ACVP2}, and motivated by further connections to the literature, we also have a worth mentioning direction, which generalizes the notion of of uniformly attainable vectors.  
    
Let $\psi\colon X\to\R$ be a continuous function, which we regard as a reference potential. For $R>0$, $N\ge1$ and $q\in\R^{N}$ with $\|q\|_\infty\le R$, let $\nu^{\psi}_{q}$ denote the unique equilibrium state of $\psi+\sum_{j\le N}q_j\phi_j$ (when it exists), and set
\begin{align*}
\cW^{\psi}_{N}(R)&:=\bigg\{\bigg(\int_{X}\phi_{1}\,d\nu^{\psi}_{q},\dots,\int_{X}\phi_{N}\,d\nu^{\psi}_{q}\bigg):\ \|q\|_{\infty}\le R\bigg\},\\
\cW^{\psi}(R)&:=\bigl\{\alpha\in\R^{\N}:\ \alpha^{(N)}\in\cW^{\psi}_{N}(R)\ \text{for all }N\ge1\bigr\}.
\end{align*}
It is clear that for $\psi\equiv0$ these are the sets $\cW_N(R)$ and $\cW(R)$ of Section~\ref{CVPS}. 

Before the statement, we recall the definition of Gibbs states. Let $T$ be a continuous map on a compact metric space $X$. We say that a measure $\mu$ on $X$ (not necessarily $T-$invariant) is a \emph{Gibbs measure} or a \emph{Gibbs state} for a function $f\colon X \to \R$ if for any sufficiently small $\epsilon > 0$ there exists a constant $K(\epsilon) \ge 1$ such that 
\begin{equation}\label{GS}
K(\epsilon)^{-1} \le \frac{\mu(B_N(x,\epsilon))}{\exp[-N P_{\rm{top}}(f) + S_Nf(x)]}\le K(\epsilon)
\end{equation}
for all $x \in X$ and $N \in \N$, where $B_N(x,\epsilon)$ is the usual Bowen ball with respect to $T$, that is, 
\[
B_{N}(x,\epsilon) = \{y \in X: d(T^k(x),T^k(y))<\epsilon \textrm{ for $0 \le k \le N-1$}\},
\]
and $P_{\rm{top}}(f)$ is the topological pressure of $f$ with respect to $T$. 

The following result generalizes Theorem G in \cite{GKLM18} for countably many constraints inside the set of invariant measures. It can be established following the proofs of Theorem \ref{ACVP2} and Proposition \ref{prop:ex1}, with the obvious modifications now taking into consideration uniformly attainable vectors relative to the fixed potential $\psi$, together with the suitable regularity assumptions on the observable.  

\begin{corollary}\label{GIBBS}
Let $T$ be a subshift of finite type, an Anosov diffeomorphism, or a uniformly expanding map of a compact metric space $X$, which is also topologically mixing. Suppose $\Phi\colon X\to\R^{\N}$ is an observable with Lipschitz coordinate functions $\phi_j$ satisfying $\sum_{j\ge1}\|\phi_j\|_{\infty}<\infty$ and $\sum_{j\ge1}\rm{Lip}(\phi_j)<\infty$, and also let $\psi\colon X\to\R$ be Lipschitz. If the vector $0 \in \cW^{\psi}(R)$ for some $R > 0$, then the measure maximizing $h_\mu(T)+\int_X\psi\,d\mu$ over $\cM(0)$ is a Gibbs state for an observable of the form 
\[
\psi+\sum_{j\ge1}q^*_j\phi_j \quad \textrm{for some $q^{*} = (q_1^{*}, q_2^{*}, \dots)$ with $\|q^{*}\|_{\infty} \le R$}. 
\]
\end{corollary}

\end{document}